\documentclass[pre,doublespace]{revtex4-1}
\usepackage[english]{babel}
\usepackage{amssymb,amsmath,graphicx,amsthm,wrapfig,caption,subcaption,epstopdf}
\usepackage{color,setspace,esint}
\usepackage{graphicx}
\usepackage{mathtools}
\usepackage{natbib}
\usepackage{hyperref}

\theoremstyle{plain}

\newcommand{\eps}{\epsilon}
\newcommand{\bra}[1]{\left(#1\right)}
\newcommand{\Bra}[1]{\left[#1\right]}
\newcommand{\BRA}[1]{\left\{#1\right\}}

\usepackage{fancyhdr}
\begin{document}

\title{Spatially localized self--assembly driven by electrically charged phase separation}

\author{Nir Gavish}
\affiliation{Department of Mathematics, Technion--Israel Institute of Technology, Technion City 32000, Israel}

\author{Idan Versano}
\affiliation{Department of Mathematics, Technion--Israel Institute of Technology, Technion City 32000, Israel}

\author{Arik Yochelis}
\affiliation{Department of Solar Energy and Environmental Physics, Swiss Institute for Dryland Environmental and Energy Research, Blaustein Institutes for Desert Research (BIDR), Ben-Gurion University of the Negev, Sede Boqer Campus, 8499000 Midreshet Ben-Gurion, Israel}

\date{\today}

\begin{abstract}
Self--assembly driven by phase separation coupled to Coulombic interactions is fundamental to a wide range of applications, examples of which include soft matter lithography via di-block copolymers, membrane design using polyelectrolytes, and renewable energy applications based on complex nano-materials, such as ionic liquids. The most common mean field framework for these problems is the non-local Cahn-Hilliard (a.k.a. Ohta-Kawasaki) framework. In this work, we study the emergence of spatially localized states in both the classical and the extended Ohta-Kawasaki model. The latter also accounts for: (\textit{i}) asymmetries in long-range Coulomb interactions that are manifested by differences in the dielectric response, and (\textit{ii}) asymmetric short-range interactions that correspond to differences in the chemical potential between two materials phases. It is shown that in one space dimension (1D) there is a multiplicity of coexisting localized solutions, which organize in the homoclinic snaking structure, bearing similarity to dissipative systems. In addition, an analysis of 2D extension is performed and distinct instability mechanisms (related to extended and localized modes) of localized stripes are discussed with respect to model parameters and domain size. Finally, implications to localized hexagonal patterns are also made. The insights provide an efficient mechanistic framework to design and control localized self--assembly that might be a plausible strategy for low cost of nano electronic applications, i.e., a rather simple nano scale fabrication of isolated morphologies.    
\end{abstract}
	
\maketitle
	
\section{Introduction}\label{sec:intro}

Self--assembly is prominent phenomena that is frequently being exploited to design patterns in variety of soft matter applications at a wide range of scales, examples of which include active layers in organic photovoltaics~\cite{ChFi:2010,OPV_Rev:2011}, membranes in fuel cells~\cite{promislow2009pem,mauritz2004state}, and multi-component fluids, such as surfactant or polymer based emulsions and micro-emulsions~\cite{tsori2009polymers}. Unlike self--organization in systems that are driven far from thermal equilibrium, such as reaction diffusion media, self-assembly models are variational, i.e., systems that are characterized by energy decrease over time. The phenomenology of self--assembly is well described by two prototypical mean field models: the Cahn--Hilliard theory for phase separation and crystal phase field for formation patterns with a characteristic spatial scale~\cite{chen2002phase,golovin2003self}. However, both frameworks rely on local interactions and neglect the impact of long-range Coulombic forces.  The latter, however, are inherently present in a wide range of chemical mixtures, examples of which include polyelectrolytes, inverse micelles and colloids, di-block copolymers, and ionic liquids, see~\cite{yochelis2015coupling} and the references therein. 

Motivated by self--assembly in polymeric compositions Ohta and Kawasaki (OK) have suggested a framework that incorporates both the short-- and the long--range interactions~\cite{ohta1986equilibrium}; the OK equation is also often referred to as the non--local Cahn--Hilliard approach. To date, not only that the OK model is being frequently employed for studies of polymers~\cite{choksi2009phase,tsori2009polymers}, but also has been recently suggested as basis for a general theory for morphology development driven by competing short-range, and long-range nonlocal interactions \cite{gavish2016theory}. In the current form, however, the model does not account for symmetry breaking properties of the coexisting phases, such as asymmetry in the double well potential and permittivity. In what follows, we study an extended version of the OK model and specifically investigate the emergence of spatially localized solutions, having properties as reviewed in~\cite{knobloch2015spatial} and the references therein.

While in the OK context most studies and applications focused on the properties of spatially periodic patterns and have described a rich variety of patterns~\cite{choksi20112d,choksi2006periodic}, spatially localized patterns such as isolated stripes (e.g.t, wires), are significant for integrated circuits as nano scales~\cite{park2003enabling,stoykovich2007directed,kim2009block,farrell2010self}. To date, such wires are being generated via polymer lithography, however, not only that the latter methods are cumbersome and expansive, stability is also a notorious issue~\cite{stoykovich2007directed}. It would be thus plausible to induce a stable localized self--assembly from a single perturbation as an alternative to multi--step lithography. However, theory of localized states in the context of self--assembly is missing although these attract much of interest in applications, such as nonlinear optics, self--organization of quantum dots, convection, and reaction--diffusion systems, and analysis due to the intriguing phenomena of homoclinic snaking~\cite{knobloch2015spatial}. Specifically, while many insights have been provided via models with local interactions (e.g., the Swift--Hohenberg model~\cite{knobloch2015spatial}), the existence and stability properties of localized states under Coulombic interactions are in general missing~\cite{glasner2010spatially}.

Here we aim to investigate the emergence of localized states in both the classical and extended Ohta--Kawasaki model and provide a basic understanding of their formation mechanism. The focus is on homoclinic snaking phenomena in 1D and transverse instabilities of localized stripes with a limiting discussion also on localized hexagonal stripes. The paper is structured as following: First, in section~\ref{sec:EOK}, we describe the Ohta--Kawasaki framework, then in Section~\ref{sec:anal}, we perform linear and weakly nonlinear analyses (deriving the Newell--Whitehead--Segel equation) and show the presence of super-- and sub--critical bifurcations. Thirdly, we study in Section~\ref{sec:Loc} the properties of spatially localized states in 1D and 2D, and finally discuss the results by focusing on the not-slanted snaking properties as could have been expected from analyses of conserved systems without Coulombic interactions~\cite{dawes2008localized,thiele2013localized}. 

\section{Extended Ohta–Kawasaki framework}\label{sec:EOK}

\begin{figure}[tp]
	\includegraphics[width=0.45\textwidth]{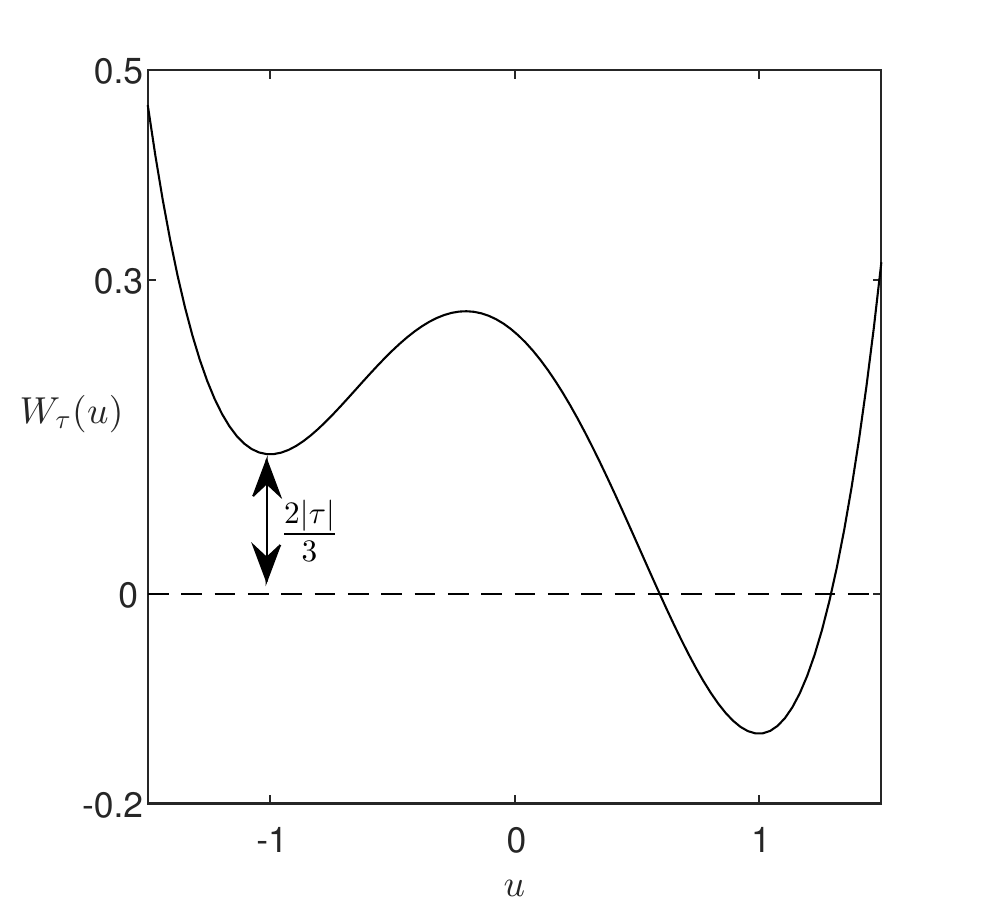}
	\caption{An example for a tilted double-well potential computed from~\eqref{eq:W_tau(u)} with $\tau=-0.2$.}\label{fig:w_tau}
\end{figure}

The competition between short-- and long--ranged Coulomb forces can be introduced by considering a mixture of two immiscible and counter-charged phases that incorporate repelling and attractive interactions, respectively. In addition to the basic OK formulation, we allow phase properties also to differ, e.g., their permittivities and/or free energies. Letting $\Omega$ be a bounded domain, the extended Ohta-Kawasaki (EOK) functional~\cite{orizaga2016instability} reads as
	\begin{equation}\label{eq:energy functional}
	\mathcal{F}=\underbrace{E_0\int_{\Omega}\frac{\varepsilon^2}{2}|\nabla u|^2+W_{\tau}(u) \ \mathrm{d}\Omega
	}_{E_{CH}}+\underbrace{\int_{\Omega}q(p(u)-n(u))\phi -\frac{1}{2}\epsilon_0(u)|\nabla \phi|^2 \ \mathrm{d}\Omega}_{E_{C}},
	\end{equation}
where $E_{CH}$ is the Cahn-Hilliard free energy associated with phase separation and $E_{C}$ stands for the long-range Coulomb forces with electrical potential $\phi$ that satisfies the Poisson's equation
	\begin{equation}\label{eq:poisson}
	-\nabla \cdot \Bra{\epsilon_0(u)\nabla \phi})=q\Bra{p(u)-n(u)}.
	\end{equation}
Furthermore, in~\eqref{eq:energy functional}, the distinct phases are given by the order parameter $u$ which is associated with the relative molar density of the positively ($u=1$) and negatively ($u=-1$) charged domains, $\eps_0(u)$ is the permittivity, $p(u):=c_p\frac{1+u}{2}$ and $n(u):=c_n\frac{1-u}{2}$ are total volume densities of the positively and negatively charged phases, respectively (with $c_p,c_n>0$), $q$ is the electrical charge, $\varepsilon$ represents the characteristic length of the interface generated by the phase separation, $E_0$ is a characteristic interfacial energy between the phases, and $W_{\tau}(u)$ is the tilted double-well potential (Figure \ref{fig:w_tau}):
\begin{equation}\label{eq:W_tau(u)}
W_{\tau}(u):=\frac{(1-u^2)^2}{4}+\tau\frac{3u-u^3}{3}.
\end{equation}
For simplicity, we assume a relatively weak difference in permittivity, i.e., linear order in $u$~\cite{orizaga2016instability},
{ 
\[
\epsilon_0(u)=Au+B,
\]
}and employ Neumann boundary conditions on $\partial \Omega$
\begin{equation}
\frac{\partial u}{\partial n}=\frac{\partial \phi}{\partial n}=\frac{\partial\, \nabla^2 u}{\partial n}=0,
\end{equation}
where $n$ is the outer unit normal to $\Omega$. These boundary conditions imply mass conservation and global electroneutrality
\begin{equation}\label{eq:mass_constraint}
\frac{\mathrm{d}}{\mathrm{d}t}\int_{\Omega}u \ \mathrm{d}\Omega =0,\quad  \quad \int_{\Omega}p(u)-n(u) \ \mathrm{d}\Omega=0,
\end{equation}
respectively. Using the notations for $p$ and $n$ we then retrieve the total mass
\begin{equation}\label{eq:mass_def}
m:=\frac{1}{|\Omega|}\int_{\Omega}u \ d\Omega=\frac{c_n-c_p}{c_n+c_p}.
\end{equation}
Finally, by combining~\eqref{eq:poisson} with $H^{-1}$ gradient-flow approach for the phase evolution~\cite{glasner2010spatially},
\begin{equation}
\frac{\mathrm{d}}{\mathrm{d}t}E[u(t)]=-D\int_{\Omega}\left |\nabla \frac{\delta E}{\delta u} \right |^2\leq 0,
\end{equation}
where $D>0$ is the diffusion constant, and by introducing rescaled variables $\tilde u=u+m,$
$$\tilde x=\frac{x}{\lambda_D}, \qquad \tilde t=\frac{DE_0}{\lambda_D^2}t, \qquad \tilde\phi=\frac{q}{K_B T}\phi,  \quad \tilde \epsilon_0 =\frac{\epsilon_0}{B},\qquad \lambda_D=\sqrt{\frac{\epsilon_0(0) K_BT}{q^2\bar{c}}},\quad \bar{c}=\frac{c_p+c_n}2,$$
we obtain the dimensionless form of the EOK equation (presented after omitting the tildes),
\begin{subequations}\label{eq:intro_pde1}
	\begin{equation}\label{eq:intro_pde_a}
	\frac{\partial u}{\partial t}=\nabla^2 \Bra{(3m-\tau)u^2+u^3-(1-3m^2+2\tau m)u-\gamma \nabla^2 u+\sigma \phi},
	\end{equation}
	\begin{equation}\label{eq:intro_pde_b}
	-\nabla \cdot \BRA{\Bra{a(u+m)+1}\nabla \phi}=u,
	\end{equation}
 	where the new parameters are given by 
\end{subequations}
	\[
	\sqrt{\gamma}=\frac{\varepsilon}{\lambda_D},\quad \sigma= \frac{\bar{c}k_BT}{E_0},\quad a=\frac{A}{B}.
	\]
In what follows, we use~\eqref{eq:intro_pde1} to study spatially localized states and chose $\gamma$ as a control parameter, since it expresses the balance between the short- and long-range forces. Namely, a balance between the characteristic length of the interfacial energy, $\varepsilon$, and the Coulombic (a.k.a. Debye) screening length $\lambda_D$. Nevertheless, the results are not limited to $\gamma$ and can be generalized to any other parameter.

\section{Linear and weakly nonlinear theory}\label{sec:anal}

Equations~\eqref{eq:intro_pde1} has uniform \textit{trivial} solutions and we perform here, a linear analysis in both time and space to identify the bifurcation onsets that are essential for the symmetry breaking and perform a weakly nonlinear analysis to capture the parameter space in which localized states can be anticipated.

\subsection{Temporal and spatial linear analysis}\label{sec:lin_anal}

\begin{figure}[tp]
	\includegraphics[width=0.5\textwidth]{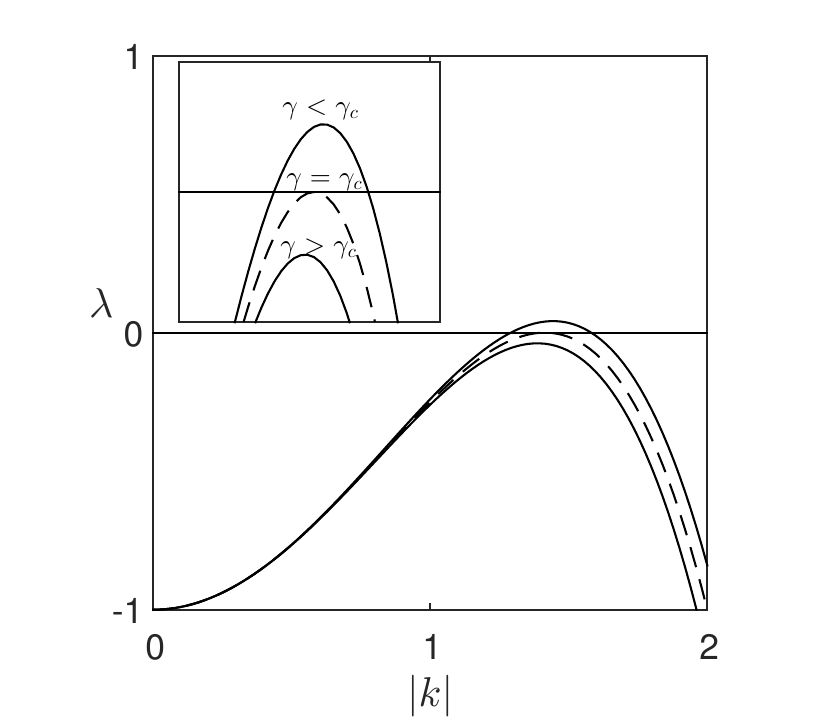}
	\caption{Dispersion relation~\eqref{eq:dispersion1} computed from~\eqref{eq:intro_pde1} for $a=\tau=0$. The inset zooms into the finite wavenumber region $k=k_c$ at ($\gamma=\gamma_c$, dashed line), above ($\gamma<\gamma_c$) and below ($\gamma>\gamma_c$) the instability onset, where $\gamma_c$ and $k_c$ are given in~\eqref{eq:gamma_c} and~\eqref{eq:k_c}, respectively.}\label{fig:disper}
\end{figure}

Let us consider solutions of~\eqref{eq:intro_pde1} on an infinite domain. Finite domain effects will be discussed in a subsequent section. At the linear order about the uniform trivial state, the solutions to~\eqref{eq:intro_pde1} are approximated as
\begin{equation}\label{eq:lin_u}
	u \simeq e^{\lambda t+ikx}+\mathrm{complex\, \, conjugated},
\end{equation}
with
\begin{equation}\label{eq:lin_phi}
\phi \simeq \frac{u}{(am+1)k^2}+\mbox{const},
\end{equation}
where $\lambda$ is the growth rate of perturbations associated with the wavenumber $k$. Substituting~\eqref{eq:lin_u} and~\eqref{eq:lin_phi} into~\eqref{eq:intro_pde1}, we obtain the dispersion relation
\begin{equation}\label{eq:dispersion1}
\lambda=(1-3m^2+2\tau m)k^2-\gamma k^4-\frac{\sigma}{am+1}.
\end{equation}
The instability onset is associated with a critical parameter for which $\lambda(k=k_c)=0$, $\lambda(k\ne k_c)<0$ and $\mathrm{d} \lambda/ \mathrm{d}k=0$ at $k=k_c$. In what follows, we use $\gamma$ as a control parameter and study the parameter space spanned by $(a,\tau)$; for simplicity we use in numerical computations $\sigma=1$. Consequently, we identify a finite wavenumber instability associated with
\begin{eqnarray}
\label{eq:gamma_c}
\gamma_c&=&\frac{(am+1)\eta^2}{4\sigma}, \\ 
k_c&=&\sqrt{\frac{\eta}{2\gamma_c}},
\label{eq:k_c}
\end{eqnarray}
and
\begin{equation}
\eta:=1-3m^2+2\tau m>0.
\end{equation}
Figure~\ref{fig:disper}, shows the dispersion relation at the instability onset (solid line) and additional two curves below and above the onset, where latter ($\gamma<\gamma_c$) is identified with a finite band of wavenumbers around $k=k_c$, for which $\lambda>0$. 

In fact, the finite wavenumber bifurcation has been already identified in the OK system ($a=\tau=0$) by Shiwa~\cite{shiwa1997amplitude} and from the above analysis the qualitative results persist also with $a,\tau \neq 0$. However, since our interest is in localized states, we examine also the configuration \textit{spatial} eigenvalues which is a necessary condition for exponentially decaying solutions as $x \to \pm \infty$~\cite{champneys1998homoclinic}, i.e., the existence of a Hamiltonian Hopf bifurcation in space.
Thus, to analyze the spatial dynamics, we use the time-independent version of \eqref{eq:intro_pde1} which forms six spatial first order ordinary differential equations:
\begin{subequations}\label{eq:dynamical system 6}
\begin{eqnarray}
	u' &=& v, \\
	v' &=& w, \\
	w' &=& z, \\
	\gamma z' &=& 3(2uv^2+u^2 w)+(6m-2\tau)(v^2+uw)-\eta w-\sigma\frac{u+av \psi}{a\bra{u+m}+1} , \\
	\phi' &=& \psi, \\
	\psi' &=& - \frac{u+av \psi}{a\bra{u+m}+1},
\end{eqnarray}
\end{subequations}
where prime denotes derivative with respect to $x$. Linearizing again about the trivial solution, we obtain six eigenvalues, among which two remain always trivial ($s=0)$ and the other four are given by
\begin{equation}\label{eq:spatial eigenvalues}
	s^2_\pm= \frac{-\eta}{2\gamma}\pm\frac{1}{2\gamma}\sqrt{\eta^2-4\sigma\frac{\gamma}{am+1}}.
\end{equation}
For $\gamma=\gamma_c$, $\eta>0$, and $0\leq a \leq 1$ the eigenvalues have a double multiplicity at the imaginary axis:
\begin{equation}\label{eq:spatial eig_c}
s^2_c= \frac{-2\sigma}{\bra{am+1}\bra{1-3m^2+2\tau m}} <0,
\end{equation}
while in the vicinity of $\gamma_c$, $\gamma-\gamma_c \sim \epsilon\ll \gamma_c$, the real and imaginary parts in~\eqref{eq:spatial eigenvalues} are determined by  
	\begin{equation}\label{eq:Hopf}
	s_\pm \sim \sqrt{ \frac{\sigma\epsilon}{\eta(am+1)}+O(\epsilon^2)} \pm i \sqrt{\eta+\frac{\sigma\epsilon}{\eta(am+1)}+O(\epsilon^2)},
	\end{equation}
where, for simplicity we ignore the eigenvalue multiplicity. Following~\eqref{eq:Hopf}, we indeed obtain at $\gamma=\gamma_c$ the Hamiltonian--Hopf bifurcation, where for $\epsilon<0$ ($\gamma<\gamma_c$)
the eigenvalues split on the imaginary axis and for $\epsilon>0$ ($\gamma>\gamma_c$) the eigenvalues upon splitting become complex and form a quartet in the real-imaginary plane~\cite{champneys1998homoclinic}. The latter case, designates also temporal stability of the uniform solution and implies a sub--critical bifurcation to localized states, bearing similarity to dissipative systems even in the presence of a double multiplicity of zero eigenvalues which is a characteristic of conservative systems~\cite{knobloch2015spatial}.

\subsection{Weakly nonlinear theory: Newell-Whitehead-Segel amplitude equation}\label{sec:NWS}

Following the possible existence of localized states as sub--critical solutions, we perform here a weakly nonlinear analysis to study their emergence and as a generalization of the OK model, also for periodic solutions that bifurcate in the super--critical direction.
Using a standard multiple time scale method~\cite{hoyle2006pattern} and the ansatz
\begin{eqnarray}
u &\sim& \sqrt{\epsilon}  A(X,Y,T)e^{ik_c x}+c.c.+ O\bra{\eps}  \\
\phi &\sim & \sqrt{\epsilon}Z(X,Y)+ \left [\sqrt{\epsilon} B(X,Y)e^{ik_c x}+c.c \right]+ O\bra{\eps} ,
\end{eqnarray}
we obtain a Newell-Whitehead-Segel (NWS) amplitude equation (see appendix for details)

\begin{equation}\label{eq:WS}
\frac{\partial A}{\partial T}=\hat{\gamma}k_c^4A+f|A|^2A+\left(\sqrt{2\eta}\frac{\partial}{\partial X}-i\sqrt{\gamma_c}\frac{\partial^2}{\partial Y^2} \right)^2A.
\end{equation}
where $A,B$ are complex amplitudes, $X=\sqrt{\epsilon}x$, $Y=\epsilon^{\frac{1}{4}}y$, $T=\epsilon t$, $\epsilon = |\gamma-\gamma_c|\ll 1$, $\hat{\gamma}:={\rm sgn}(\gamma_c-\gamma)$, $c.c.$ stand for complex conjugate, and 
\begin{equation}\label{eq:f}
\begin{split}
f=&\frac{2\sigma}{9\eta^2(am+1)^3} (198a^2m^4-108a^2m^3\tau-3a^2m^2-18a^2m\tau+468am^3-318am^2\tau+36am\tau^2-5a^2 \\
& -60am+2a\tau+225m^2-150m\tau+16\tau^2-27).
\end{split}
\end{equation}
Transforming back to fast variables, we obtain:
\begin{equation}\label{eq:WS_fast_variablesd}
	\frac{\partial A}{\partial t}=(\gamma_c-\gamma)k_c^4A+\frac{\gamma_c-\gamma}{{\rm sgn}(\gamma_c-\gamma)}f|A|^2A+\left( \sqrt{2\eta}\frac{\partial}{\partial x}-i\sqrt{\gamma_c} \frac{\partial^2}{\partial y^2} \right)^2 A.
\end{equation}

As expected, Eqs.~\ref{eq:WS} (and respectively~\ref{eq:WS_fast_variablesd}) gives rise to three type of solutions:
\begin{itemize}
\item Constant
\begin{equation}\label{eq:A_c}
A_p =\sqrt{ \frac{\hat{\gamma}k_c^4}{-f} }e^{i \psi}, \qquad \psi\in \mathbb{R},
\end{equation}
which corresponds periodic solutions in EOK
\begin{equation}\label{eq:A_c}
u_p = \sqrt{\frac{\gamma_c-\gamma}{-f}}k_c^2e^{i \psi} +c.c.\, .
\end{equation} 

\item Periodic 
\begin{equation}\label{eq:A_Q}
A_Q = \sqrt{\frac{2 \eta Q^2-\hat{\gamma} k_c^4}{f}}e^{iQX+i\psi}.
\end{equation}
which corresponds modulations of the periodic solutions in EOK
\begin{equation}\label{eq:u_Q}
u_Q = \sqrt{|\gamma_c-\gamma|}\sqrt{\frac{2\eta Q^2-\hat{\gamma} k_c^4}{f}}e^{i(k_c+\sqrt{|\gamma_c-\gamma|}Q)x+i\psi}+c.c.\, .
\end{equation}

\item Spatially localized
	\begin{equation}\label{eq:sech}
	A_L:=\sqrt{\frac{-2\hat{\gamma} k^4_c}{f}} {\rm sech} \bra{\sqrt{\frac{-\hat{\gamma}k_c^4}{2 \eta}} X}e^{i\psi}, 
	\end{equation}
which corresponds to solutions that decay exponentially as $x \to \pm \infty$
\begin{equation}\label{eq:u_L}
u_L=\sqrt{\gamma-\gamma_c}\sqrt{\frac{2k_c^4}{f}} {\rm sech} \left(\sqrt{\gamma-\gamma_c}\sqrt{\frac{k_c^4}{2 \eta}}x \right) e^{ik_c x+i\psi} + c.c.\,.
\end{equation}
\end{itemize}
For $\gamma<\gamma_c$ and $f<0$, branches of the spatially periodic solutions $u_p$ and $u_Q$ bifurcate in a super--critical fashion (towards the linearly unstable regime), while for $\gamma>\gamma_c$ and $f>0$ they bifurcate sub-critically (towards the linearly stable regime) together with spatially localized solutions $u_L$. Although in~\eqref{eq:u_L} $u_L$ is phase invariant, beyond-all-orders analysis shows that under even symmetry, $\psi$ eventually locks to two relative values of $0$ and $\pi$~\cite{chapman2009exponential}, as will be also confirmed numerically here. 

\subsection{Super-criticality: periodic solutions and Busse balloon}\label{sec:sup}

Before starting the analysis of spatially localized solutions, we wish to generalized the properties of periodic solutions in 1D and their respective analogue of stripes in 2D, i.e., compute the onset of secondary instabilities with respect to Eckhaus and zig--zag, and construct the Busse balloon in a similar fashion that was done for the classical OK model ($a=\tau=0$)~\cite{shiwa1997amplitude}.

\subsubsection{Eckhaus instability}\label{sec:Eckh}
\begin{figure}[tp]
	\large{(a)\includegraphics[width=0.43\textwidth]{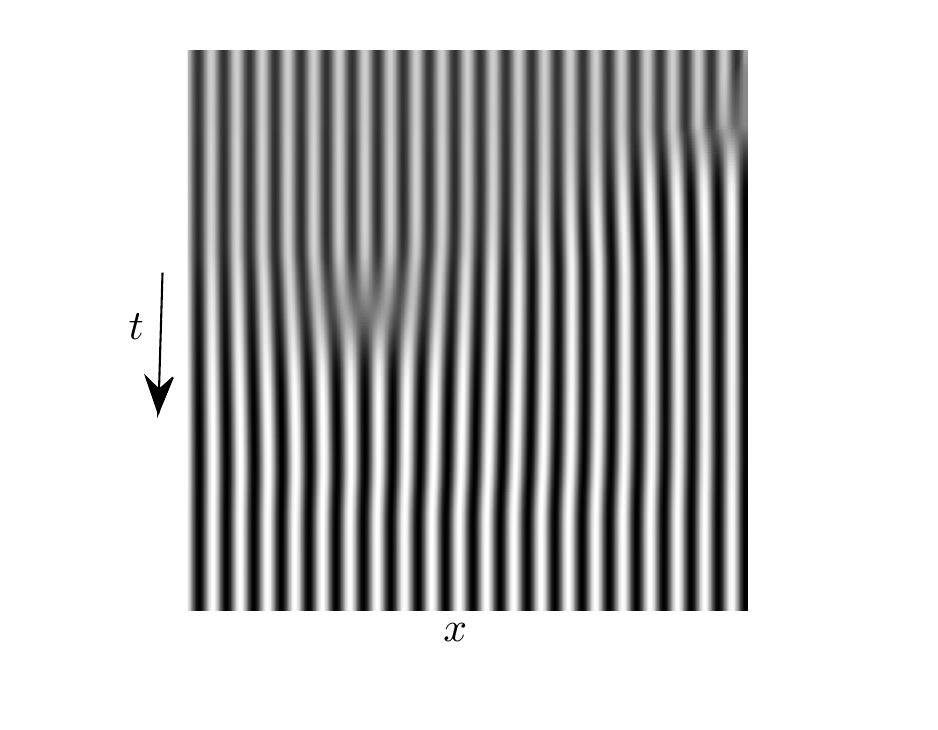}
		(b)}\includegraphics[width=0.45\textwidth]{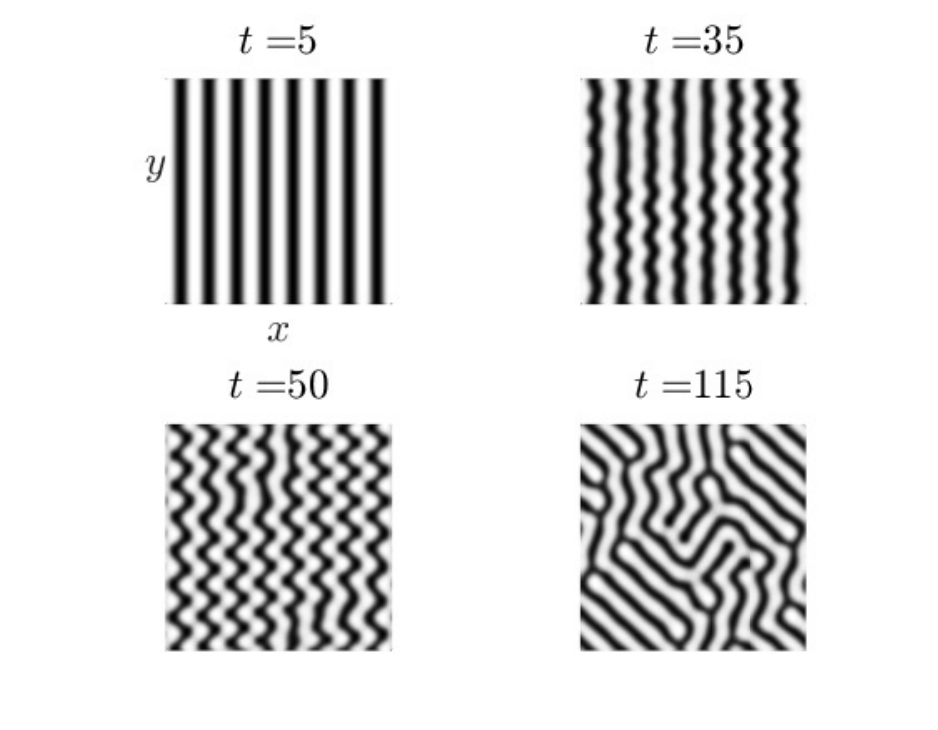}
	\caption{(a) Space--time plot showing the Eckhaus instability by direct numerical integration of~ \eqref{eq:intro_pde1}, at $\gamma=0.24$ with $\Omega=[0,20\frac{2\pi}{k_c}]$ and $t\in [0,1000]$, where $k_c\simeq 1.923$. (b) Snapshots at given times, showing the evolution of stripes to labyrinthine-type pattern due to the zigzag instability. Here $\gamma=0.1$ and $\Omega=[0,8\frac{2\pi}{k_c}]\times [0,8\frac{2\pi}{k_c}]$. Other parameters for both (a) and (b): $a=0.1$, $m=\tau=0$, and $\sigma=1$.}\label{fig:zigzag_non}
\end{figure}

The Eckhaus instability is a phase instability in $x$ that results in a weak modulation of the primary period $k_c$, so that solutions $A_Q$ are different from $A_p$ up to an $O(\sqrt{\epsilon})$ phase shift.
To compute the instability onset, we consider the ansatz
\begin{equation}\label{eq:eckhaus1}
A_{\text{pert}}:= A_Q+a_k(X)\cdot e^{\lambda T},
\end{equation}
where $a_k=F_{+k} e^{i(Q+k)X}+F_{-k} e^{i(Q-k)X}$ with $F_{\pm k} \ll 1$. Substituting \eqref{eq:eckhaus1} into \eqref{eq:WS} we obtain the linear system
\begin{equation}\label{eq:matrix1}
\left(
\begin{matrix}
-\hat{\gamma}k_c^4-2\eta (Q+k)^2+4 \eta Q^2 & -\hat{\gamma}k_c^4+2 \eta Q^2 \\
-\hat{\gamma}k_c^4+2 \eta Q^2 & -\hat{\gamma}k_c^4-2\eta (Q-k)^2+4 \eta Q^2
\end{matrix} \right )
\left(
\begin{matrix}
F_{+k} \\
F_{-k}
\end{matrix} \right )
=\lambda \left(
\begin{matrix}
F_{+k} \\
F_{-k}
\end{matrix} \right ).
\end{equation}
with eigenvalues given by
\begin{eqnarray}\label{eq:eigen2}
\lambda_0&=&2(2\eta Q^2-\hat{\gamma}k_c^4), \\
\lambda_k^{\pm}&=&2\eta Q^2-\hat{\gamma}k_c^4-2 \eta k^2 \pm \sqrt{( \hat{\gamma}k_c^4-2\eta Q^2)^2+16\eta^2 Q^2k^2}.
\end{eqnarray}
It that both $\lambda_0$ and $\lambda_{k}^{-}$ are negative, so that the Eckhaus instability corresponds to the sign of $\lambda_{k}^{+}$ which can get positive values in the range of 
\begin{equation}
Q^2<\frac{\hat{\gamma}k_c^4}{2\eta}<3Q^2.
\end{equation}
In this region the solution $u_Q$, see~\eqref{eq:u_Q}, loses stability to the modes which satisfy
\begin{equation}
|k|<\sqrt{6Q^2-\frac{\hat{\gamma}k_c^4}{\eta}},
\end{equation}
as shown by time integration in Figure~\ref{fig:zigzag_non}(a). In general, direct numerical integrations in 1D, were performed by a standard finite differences scheme where the derivatives with respect to time and space being calculated explicitly and implicitly, respectively.
In each time step, the system of finite-difference equations was solved using Newton iterations, where the initial guess was chosen to be the forward Euler approximation, see~\cite{christlieb2014high} for details.

\subsubsection{Zigzag instability}\label{sec:ZZ}
The second phase instability corresponds to a transverse direction to stripes, a.k.a. zigzag instability. Consequently, the ansatz \eqref{eq:eckhaus1} is modified as $a_l=F_{+l}e^{i(QX+lY)}+F_{-l}e^{i(QX-lY)}$. The eigen system then becomes
\begin{equation}\label{eq:matrix1}
\left(
\begin{matrix}
2\eta Q^2-2\sqrt{2\eta\gamma_c}Ql^2-\gamma_cl^4-\hat{\gamma}k_c^4 & -\hat{\gamma}k_c^4+2\eta Q^2 \\
-\hat{\gamma}k_c^4+2\eta Q^2 & 2\eta Q^2-2\sqrt{2\eta\gamma_c}Ql^2-\gamma_cl^4-\hat{\gamma}k_c^4
\end{matrix} \right )
\left(
\begin{matrix}
F_{+l} \\
F_{-l}
\end{matrix} \right )
=\lambda \left(
\begin{matrix}
F_{+l} \\
F_{-l}
\end{matrix} \right ),
\end{equation}
and the corresponded eigenvalues are
\begin{eqnarray}
\lambda^{zz}_1&=&-l^4\gamma_c-2Ql^2\sqrt{2\eta \gamma_c}, \\
\lambda^{zz}_2&=&-l^4\gamma_c-2Ql^2\sqrt{2\eta \gamma_c}+\underbrace{2(-\hat{\gamma}k_c^4+2\eta Q^2)}_{\leq 0},
\end{eqnarray}
so that $\lambda^{zz}_1>\lambda^{zz}_2$. The zigzag instability is of long wavenumber so that a stripe solution of wavenumber $k_c-\sqrt{\epsilon}Q$ will be unstable for perturbations in $y$ direction, provided that $l^2<-2\sqrt{2\eta}Q/\sqrt{\gamma_c}$. Figure \ref{fig:zigzag_non} shows the time evolution of a stripe solution with weak initial random noise. We employed for 2D direct numerical schemes the pseudo-spectral discretization following Eyre's method (Eyre used only for $a=0$ otherwise we split linear and nonlinear terms) \cite{eyre1998unconditionally}, by splitting $\mathcal{F}$ (see Eq.~\ref{eq:energy functional}) into contractive and expansive terms (exploiting that $|u|<1$).

\section{Spatially localized states}\label{sec:Loc}
Typically, bifurcating localized solutions are unstable near the bifurcation onset and gain stability via saddle node bifurcations~\cite{knobloch2015spatial}. To follow the branch of asymptotic (time independent) localized state, we employ a numerical continuation method using the package AUTO~\cite{doedel2012auto}. Unless stated otherwise we use the norm
\begin{equation}
	||u||=\sqrt{\int_{0}^{1}u^2+v^2+w^2+z^2+\phi^2+\psi^2 \ dx},
\end{equation}
to plot the bifurcations in the $(\gamma,||u||)$ plane. As has been already noted, the parameters $a$ and $\tau$ do not add any qualitative changes to the linear problem, and thus for simplicity we start with OK model and then generalize to EOK.
	
\subsection{Ohta–Kawasaki model, $a=\tau=0$}\label{sec:Loc_OK}

\subsubsection{Homoclinic snaking in 1D}\label{sec:Loc_snaking}

\begin{figure}[tp]
\includegraphics[width=\textwidth]{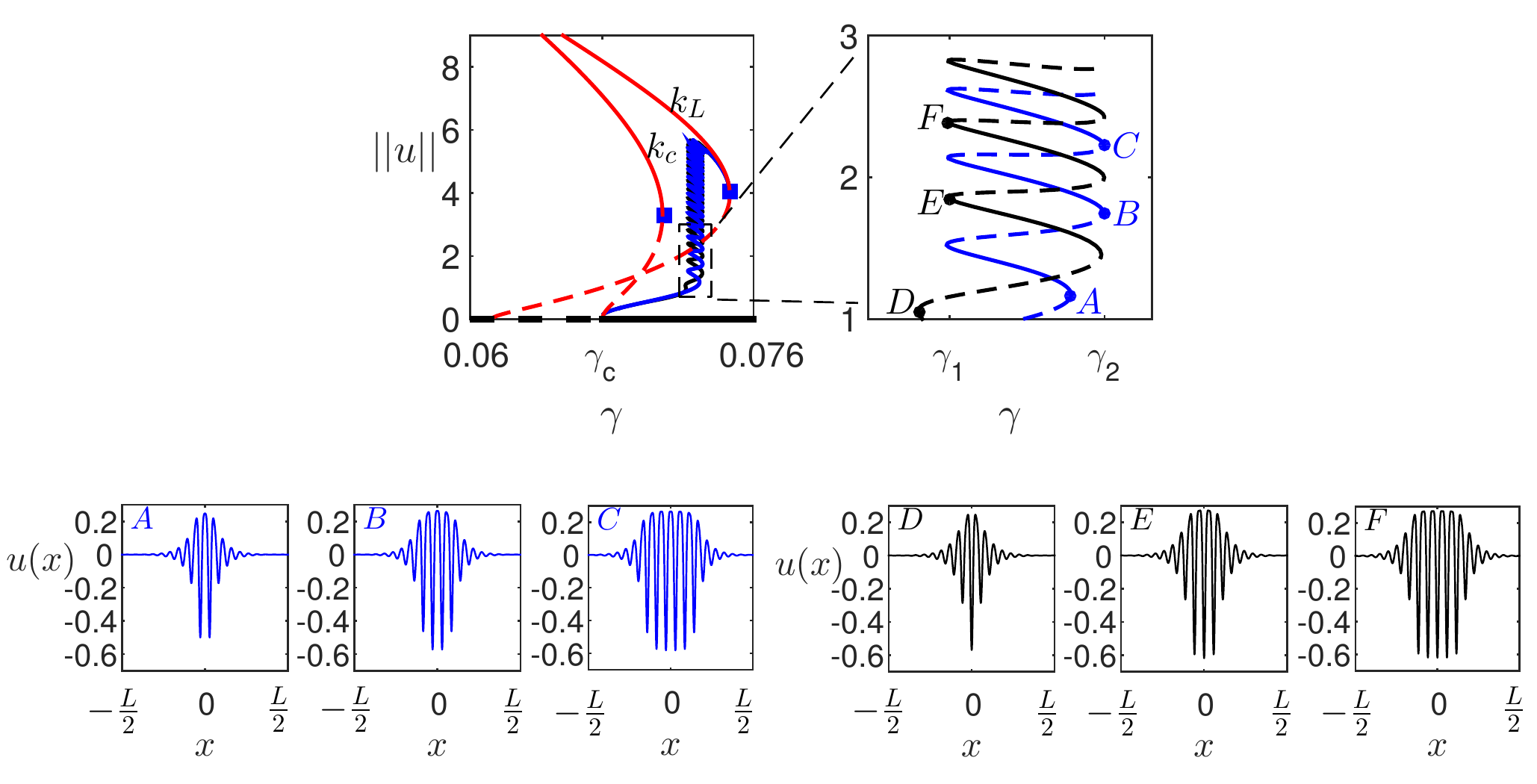}
	\caption{Top panel: Bifurcation diagram showing both periodic and localized solutions to~\eqref{eq:dynamical system 6} obtained via continuation method~\cite{doedel2012auto}. Solid and dashed lines in the zoom in, represent stable and unstable solutions, respectively. The localized states are organized in the so called homoclinic snaking structure. The primary bifurcating periodic solution is marked by $k_c$ ($P_c=\frac{2\pi}{k_c}$) while $k_L$ ($P_L=\frac{40}{35}\frac{2\pi}{k_c}$, $k_c \simeq 1.96$) stands for a periodic solution on which the branches of localized states terminate. Bottom panel: spatial profiles at respective locations as indicated as the top panel. Neumann BC have been employed for all fields except $\phi(L/2)=0$, where $L=40 \pi/k_c$. However, the latter is equivalent here to Neumann BC, as was also independently verified by an integral condition. Parameters: $m=0.4$, $a=\tau=0$, and $\sigma=1$.}\label{fig:snaking_ok}
\end{figure}

Spatially localized solutions~\eqref{eq:u_L} bifurcate in two branches toward the region where uniform state is linearly stable ($\gamma>\gamma_c$), as shown in Figure~\ref{fig:snaking_ok}. We mark these branches as $L_0$ and $L_{\pi}$ due to the relative phase difference that exists beyond all order computations~\cite{chapman2009exponential}, i.e., in \eqref{eq:u_L} $\psi=0,\pi$. These branches snake back and forth through folds, in an interval $\gamma \in [\gamma_1,\gamma_2]$. This intertwined snaking structure has a standard (vertical) form and bares similarity to models with dissipative properties, such as Swift-Hohenberg~\cite{burke2007homoclinic} and Gierer-Meinhardt~\cite{yochelis2008formation}. Within this region, each branch corresponds after an odd (even) fold to stable (unstable) localized solutions with either even or odd number of peaks while the number of peaks at each stable branch increases by two, as shown in the snaking region zoom in and the respective profiles in Figure~\ref{fig:snaking_ok}. 

The stability analysis of localized solutions is performed numerically by expanding
\begin{equation}\label{eq:1D_ansatz}
\left (
\begin{matrix}
u \\
\phi
\end{matrix}
\right)=
\left (
\begin{matrix}
u_L(x) \\
\phi_L(x)
\end{matrix}
\right)+
\left (
\begin{matrix}
\tilde{u}(x) \\
\tilde{\phi}(x)
\end{matrix}
\right)e^{\beta t},
\end{equation}
and solving at the leading order ($\tilde{u},\tilde{\phi}\ll 1$) the eigenvalue problem:
\begin{subequations}
\begin{eqnarray}
	\label{eq:stab_1D_u}
	\beta \tilde{u}&=&\sigma \tilde{\phi}''-\gamma \tilde{u}'''' -\eta \tilde{u}''+\bra{6m-2\tau}\bra{u_L\tilde{u}}''+3\bra{u_L^2\tilde{u}}'',  \\
	0&=&\bra{am+1}\tilde{\phi}''+\tilde{u}+a \bra{u_L \tilde{\phi}'}'+a\bra{\phi_L^\prime \tilde{u}}',
	\label{eq:stab_1D}
\end{eqnarray}
\end{subequations}
where primes stand for derivatives with respect to the argument. Equation~\eqref{eq:stab_1D} is solved first, using $\mu(x)=a(u_L+m)+1\ne 0$, so that
\begin{equation}
	(\mu\tilde{\phi}')'=-\tilde{u}-(a\tilde{u}\phi_L')'.
\end{equation}
To compute the stability of $(u_L,\phi_L)$ with Neumann BC, we proceed with
\begin{equation}
	\mu\tilde{\phi}'+a\tilde{u}\phi_L'=-\int_{-L}^{x}\tilde{u}(s)ds,
\end{equation}
which leads to
\begin{equation}
	\tilde{\phi}'=-\frac{1}{\mu}\left(  a\tilde{u}\phi_L' +\int_{-L}^{x}\tilde{u}(s) ds \right),
\end{equation}
and finally
\begin{equation}\label{eq:stab_phi}
	\tilde{\phi}''=\frac{au_L'}{\mu^2}\left(  a\tilde{u}\phi_L' +\int_{-L}^{x}\tilde{u}(s) ds \right)-\frac{1}{\mu}\left(\tilde{u}+a\tilde{u}'\phi_L'+a\tilde{u}\phi_L''\right).
\end{equation}
Substitution of~\eqref{eq:stab_phi} into \eqref{eq:stab_1D_u} allows to solve for $\beta$ and thus determine the stability. We note that this computation was confirmed by an independent method that we use in 2D stability analysis (by setting $k_y=0$) and by direct numerical integrations. 
	
On an infinite domain, the branches would emerge from $\gamma_c$ and snake back and forth without any limitations. However, on finite (albeit large) domains $L=n\frac{2\pi}{k_c} \gg P_c=\frac{2\pi}{k_c}$ ($n\in \mathbb{N}$), the localized solutions bifurcate from the periodic branch $P_c$ for $\gamma \gtrsim \gamma_c$~\cite{PhysRevE.78.046201}, a.k.a. a finite domain effect, and in fact understood using by Jacobi elliptic functions as solutions to~\eqref{eq:WS_fast_variablesd} in vicinity of the onset. The second impact of the finite domain effect is related to branch termination, which happens to be on a periodic solution after the domain gets filled by peaks~\cite{knobloch2015spatial}. However, unlike typical behavior in dissipative models~\cite{burke2007snakes}, here the branch terminates at a distinct periodic solution, to which we refer as $\frac{2\pi}{k_L}=\frac{35}{80}\frac{2 \pi}{k_c}$, as shown in Figure~\ref{fig:snaking_ok}.

In general, insights to formation and properties of the snaking region exploit the presence of conserved quantities (excluding reaction--diffusion systems)~\cite{knobloch2015spatial}, e.g., Hamiltonian or chemical potential. Thus, typically the snaking region is organized in a parameter space around a Maxwell-point (here $\gamma=\gamma_{Max}$), a point where the multiple peaks (a finite periodic state) inside the localized state have the same energy as the surrounding uniform state~\cite{PhysRevE.78.046201}. In the EOK case, there is an absence of a conserved quantity for~\eqref{eq:intro_pde1}, due to the nonlocal dependence of the function $\phi$ on $u$. However, since the system is conservative (although the conserved quantity is not depended locally on $u$), it is still possible to identify the Maxwell-point from energy considerations and stability of periodic solutions.

Notably, in models which account for short-range interactions, such as the Swift-Hohenberg~\cite{burke2007homoclinic} and the conserved Swift-Hohenberg, the conserved quantity is uniquely determined by the problem parameters and the mass of the initial condition. The latter dictates structure of the snaking region by local minimizers of the energy under the mass constraint, resulting in a slanted snaking pattern~\cite{thiele2013localized}. In contrast, here, the conserved quantity takes the form
\begin{equation}\label{eq:Lag}
	\nu=\eta u+(3m-\tau)u^2+u^3+\sigma \phi-\gamma \nabla^2 u,
\end{equation}
where~$\phi$ is determined only up to a constant (see appendix). The latter implies that the conserved quantity is not determined with the solution for mass, as in the cSH case~\cite{thiele2013localized}.

\begin{figure}[tp]
	\includegraphics[width=0.5\textwidth]{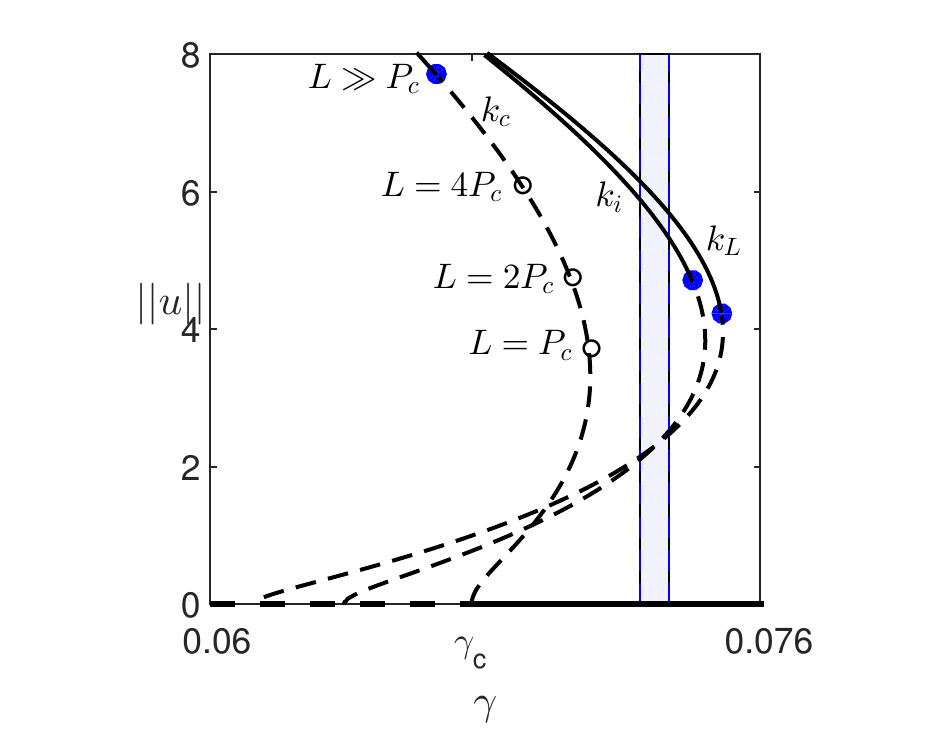}
	\caption{Branches of periodic solutions on large domains ($L$) and their relation to stability to Eckhaus; dashed lines indicate unstable solutions. The secondary two branches are related to the critical spatial period $\frac{2\pi}{k_c}$, as $\frac{2\pi}{k_i}=\frac{40}{35}\frac{2\pi}{k_c}$ and $\frac{2\pi}{k_L}=\frac{40}{36}\frac{2\pi}{k_c}$, where $k_c\simeq 1.96$. Solid circles represent the Eckhaus instability onsets while hollow circles on the $k_c$ branch correspond to the consecutive locations of the instability point when increasing the domain size. The snaking region (shaded region) is presented for reference. Parameters: $m=0.4$, $a=\tau=0$, $\sigma=1$.}\label{fig:Eckhaus}
\end{figure}

\begin{figure}[tp]
\includegraphics[width=\textwidth]{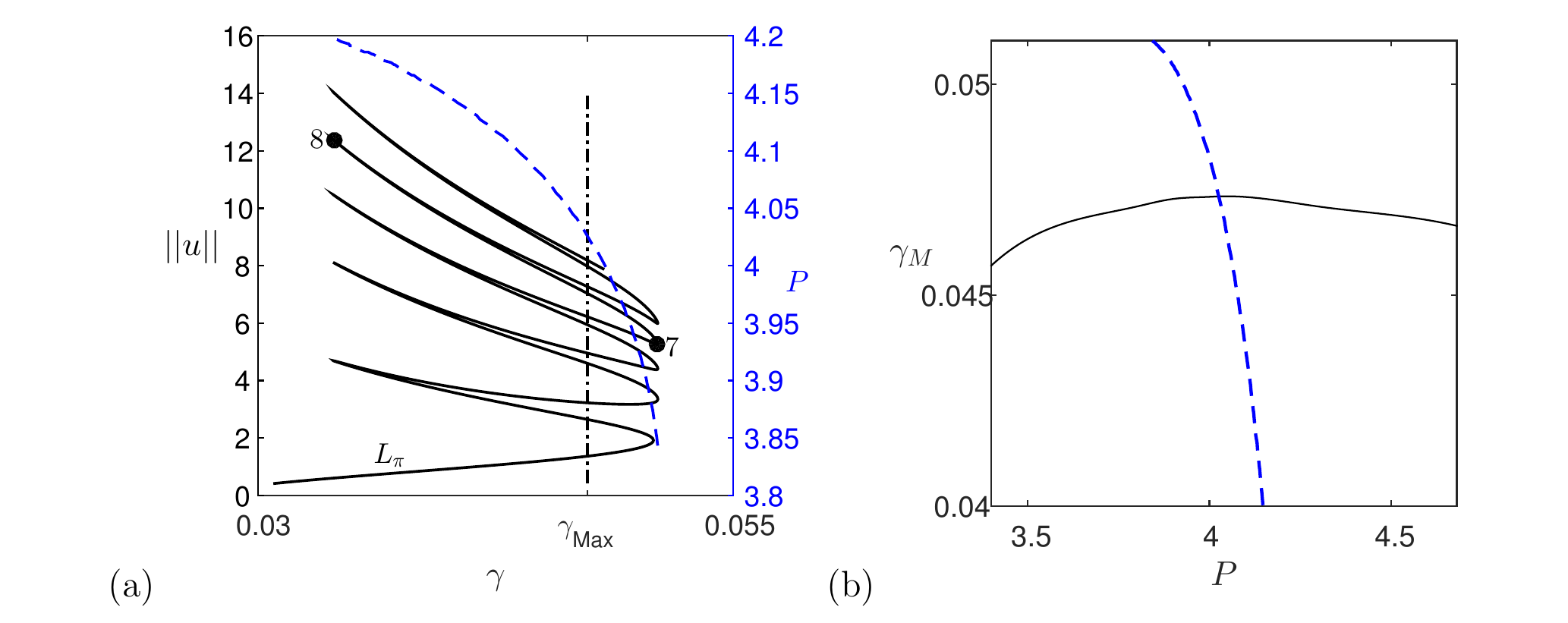}
\caption{(a) Location of the Maxwell point $\gamma_{Max}$ and the frustration effect. The latter is associated with spatial period $P$ (dashed line) variation inside periodic localized state in the snaking region, presented for example in between the seventh and eighth saddle nodes (solid line) of the localized solutions along the $L_\pi$ branch.
(b) Numerical identification of the Maxwell point, denoted by intersection of locus of cases in which periods have the same energy as the uniform state $\gamma_M$ (solid line) and the respective minimal energy of the periodic solutions (dashed line). All solutions have been obtained via the numerical continuation method. Parameters: $m=0.47$, $\tau=a=0$, and $\sigma=1$.}\label{fig:maxwell}
\end{figure}

The primary periodic solution ($k=k_c$) stabilizes at a fold where the branch direction is reversed toward decreasing $\gamma$ values. However, this solution on larger domains is Eckhaus unstable and gains ultimate stability at $\gamma<\gamma_c$ (a region in which the uniform state is linearly unstable)~\cite{yochelis2008generation,yochelis2008formation}, see solid circle on branch $k_c$ in Figure~\ref{fig:Eckhaus}. Consequently, a Maxwell point between the primary periodic state and the uniform cannot exist. The same is true for several additional secondary periods that emerge from $\gamma<\gamma_c$. Nevertheless, there exist periodic solutions that gain stability to Eckhaus above $\gamma_c$ and thus form multi--stability in the same parameter range, e.g., $k=k_i$  in Figure~\ref{fig:Eckhaus}. This coexistence has also a signature in the snaking region. Indeed, examination of profiles along the stable portion of the $L_{\pi}$ branch (in between seventh and eighth saddle nodes), shows that oscillations inside the localized state correspond to distinct periods that are selected not according to minimal energy of~\eqref{eq:intro_pde1}, see dashed line in Figure~\ref{fig:maxwell}(a). Consequently, such a coexistence designates locus of points, $\gamma=\gamma_M$, at which the periodic solutions ($P=\frac{2\pi}{k}$) has the same energy as the uniform solution, see solid line in Figure~\ref{fig:maxwell}(b). Finally, up to a numerical validity the Maxwell-point (see also $\gamma_{Max}$ in Figure~\ref{fig:maxwell}(a)), is identified by the intersection between $\gamma_M$ (solid line) and the periodic solutions (dashed line), and thus appears as a global minima in energy landscape.

\subsubsection{Localized stripes in 2D and transverse instabilities}\label{sec:Loc_2D}
	
Spatially localized states in 1D correspond to stripes in quasi-2D geometry~\cite{stoykovich2007directed} and thus for example, of interest to nanowire fabrication~\cite{farrell2010self}.	Unlike the energetic/Maxwell analogue developed above for 1D, it was shown that a similar construction for 2D and thus existence of stable localized stripes on infinite domain is unlikely to form~\cite{choksi20112d}. Nonetheless, our interest is in finite domains, for which stable localized stripes are known to persist~\cite{burke2007homoclinic,lloyd2008localized}. 
	
For the analysis, we consider stripe solutions that correspond to stable 1D localized solutions obtained via 1D continuation, as described in Section~\ref{sec:Loc_snaking}. In 2D, additional transverse instabilities often appear and thus, reduce the parameter region where the solutions gain linear stability~\cite{burke2007homoclinic,lloyd2008localized}. Linear stability of localized stripes $(u_L(x,y),\phi_L(x,y))$ can be obtained by using the anstaz
\begin{equation}\label{eq:2D_ansatz}
\left (
\begin{matrix}
u \\
\phi
\end{matrix}
\right)=
\left (
\begin{matrix}
u_L \\
\phi_L
\end{matrix}
\right)+
\left (
\begin{matrix}
\tilde{u}(x) \\
\tilde{\phi}(x)
\end{matrix}
\right)e^{\beta t+ik_y y}.
\end{equation}
Substituting~\eqref{eq:2D_ansatz} in~\eqref{eq:intro_pde1}, we obtain
\begin{subequations}\label{eq:polynomial_eigenvalue}
	\begin{eqnarray}
	\beta\tilde{u}&=&\sigma\mathcal{D}\tilde{\phi}-\gamma \mathcal{D}^2\tilde{u}-\eta \mathcal{D}\tilde{u}+\bra{6m-2\tau}\mathcal{D}u_L\tilde{u}+3\mathcal{D}u_L^2\tilde{u}, \\
	0&=&(am+1)\mathcal{D}\tilde{\phi}+\tilde{u}+a\partial_x(u_L\partial_x\tilde{\phi})+a\partial_x(\tilde{u}\partial_x\phi_L)-ak_y^2u_L\tilde{\phi},
	\end{eqnarray}
\end{subequations}
where $\mathcal{D}:=\partial_{xx}-k_y^2$. To solve this polynomial eigenvalue problem, we rewrite~\eqref{eq:polynomial_eigenvalue} as
		\begin{equation}
		\mathcal{S}\left (
		\begin{matrix}
		\tilde{u} \\
		\tilde{\phi}
		\end{matrix}
		\right) :=		
		{\left [
			\begin{matrix}
			\mathcal{L}_1(u_L,\phi_L) & \mathcal{L}_2(u_L,\phi_L) \\
			\mathcal{L}_3(u_L,\phi_L) & \mathcal{L}_4(u_L,\phi_L)
			\end{matrix}
			\right]}
		\left (
		\begin{matrix}
		\tilde{u} \\
		\tilde{\phi}
		\end{matrix}
		\right)=\beta
		\left (
		\begin{matrix}
		\tilde{u} \\
		0
		\end{matrix}
		\right),
		\end{equation}
		where $\mathcal{L}_{1,2,3,4}$ are linear operators. Applying the projection 
		\begin{equation}
		\left (
		\begin{matrix}
		\tilde{u} \\
		0
		\end{matrix}
		\right)=\mathcal{P}
		\left (
		\begin{matrix}
		\tilde{u} \\
		\tilde{\phi}
		\end{matrix}
		\right)
		\end{equation}
		reduces~\eqref{eq:polynomial_eigenvalue} to a generalized eigenvalue problem that is solvable numerically
		\begin{equation}\label{eq:eignevalue_problem_2D_stripes}
		\Bra{\mathcal{S}-\beta \mathcal{P}}
		\left (
		\begin{matrix}
		\tilde{u} \\
		\tilde{\phi}
		\end{matrix}
		\right)=
		\left (
		\begin{matrix}
		0 \\
		0
		\end{matrix}
		\right).
		\end{equation}

	\begin{figure}[tp]
		\large{
		(a)\includegraphics[width=0.45\textwidth]{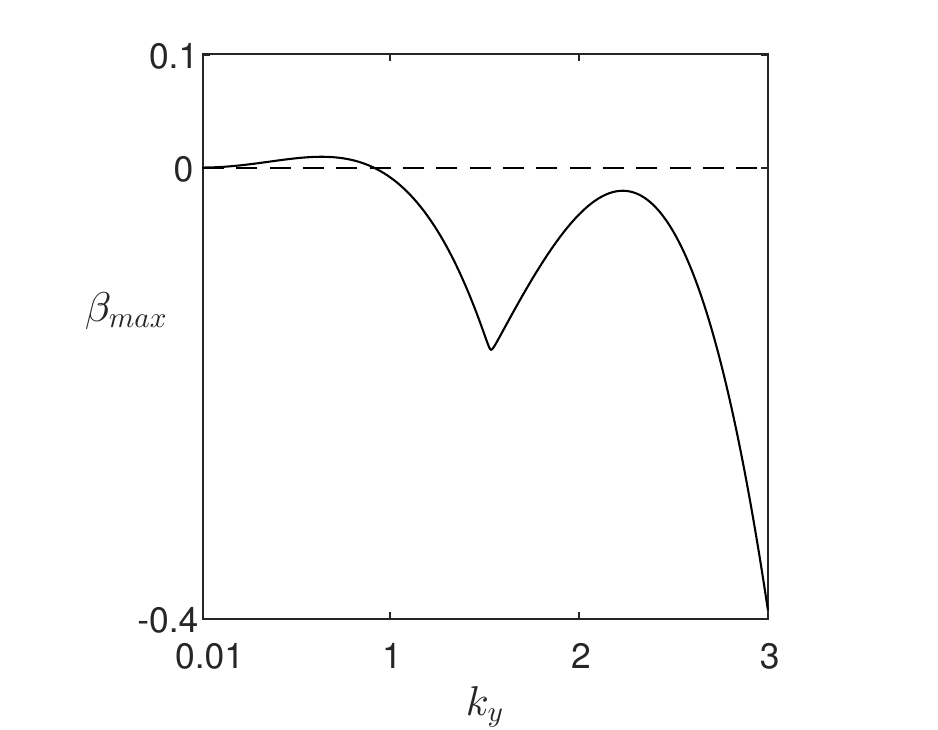}
		(b)}\includegraphics[width=0.45\textwidth]{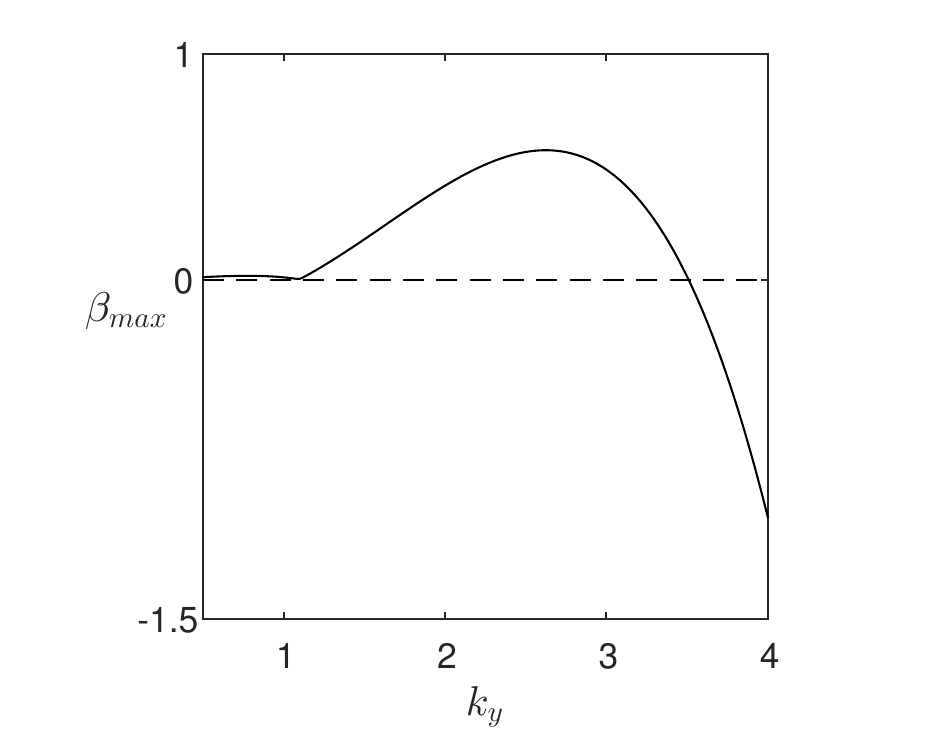}
		\caption{Dispersion relation showing the transverse instabilities of stripes obtained from~\eqref{eq:eignevalue_problem_2D_stripes}. (a) A long wave-number instability (body) at $\gamma=0.0242$ and (b) a finite wave-number instability corresponding to the wall mode at $\gamma=0.0437$. Parameters: $m=0.5$, $\tau=a=0$, and $\sigma=1$.}
		\label{fig:dispersion_wall_body}
	\end{figure}
	
\begin{figure}[tp]
	\includegraphics[width=\textwidth]{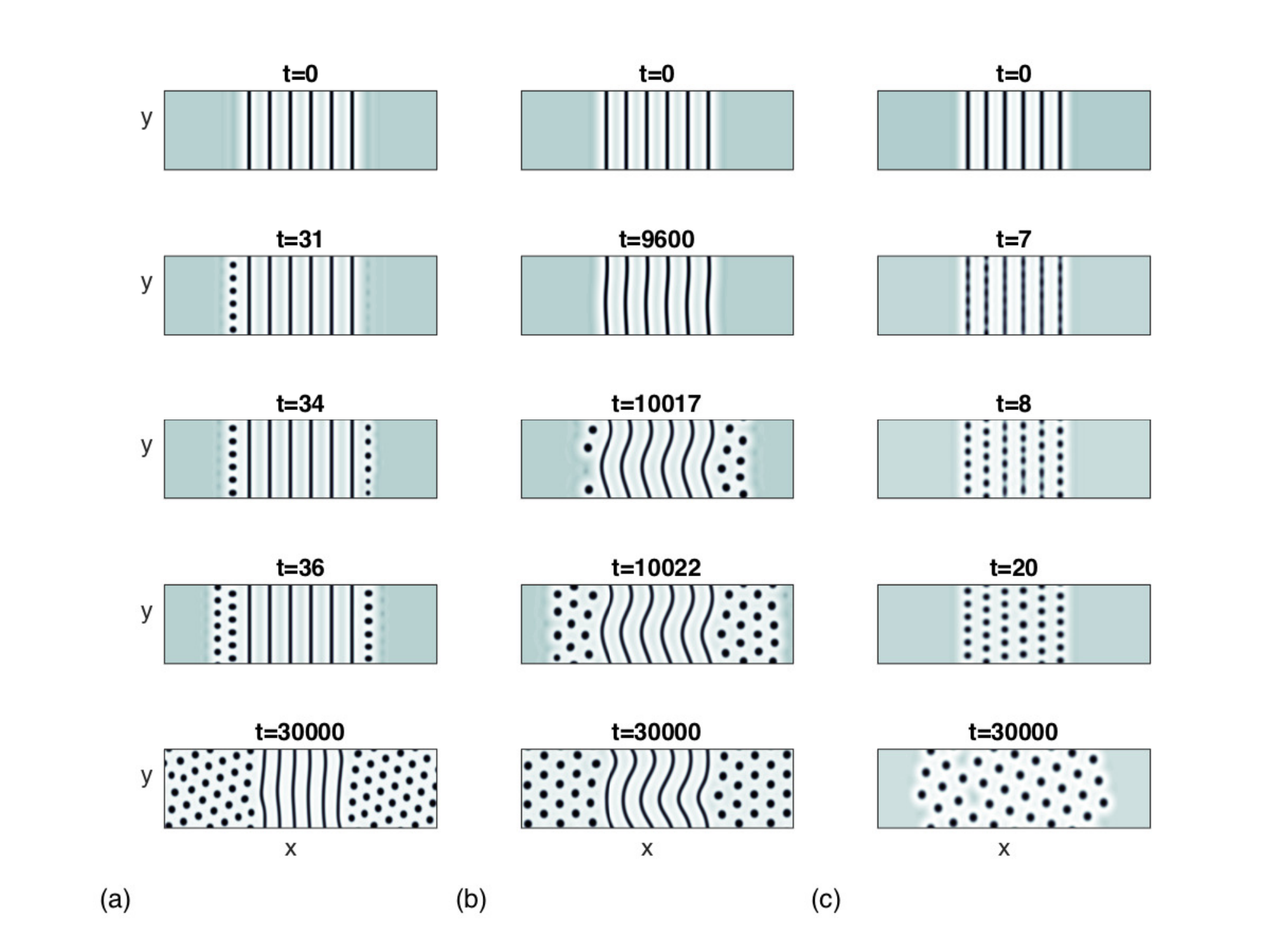}
	\caption{Snapshots at given times, showing the evolution of localized stripes due to three types of instability mechanisms: (a) body at $\gamma=0.0242$, (b) wall at $\gamma=0.0199$, and (c) localized--body at $\gamma=0.0437$.
	Other parameters: $m=0.5$, $a=\tau=0$, $\sigma=1$.  Solutions are calculated on domain~$(x,y)\in[0,L]\times[0,15]$, where~$L=80\pi/k_c$ and~$k_c\simeq 2.822$, but, for clarity, presented on the subdomain~$[L/6,5L/6]\times[0,15]$.}\label{fig:TI}
\end{figure}

\begin{figure}[tp]
	\includegraphics[width=0.7\textwidth]{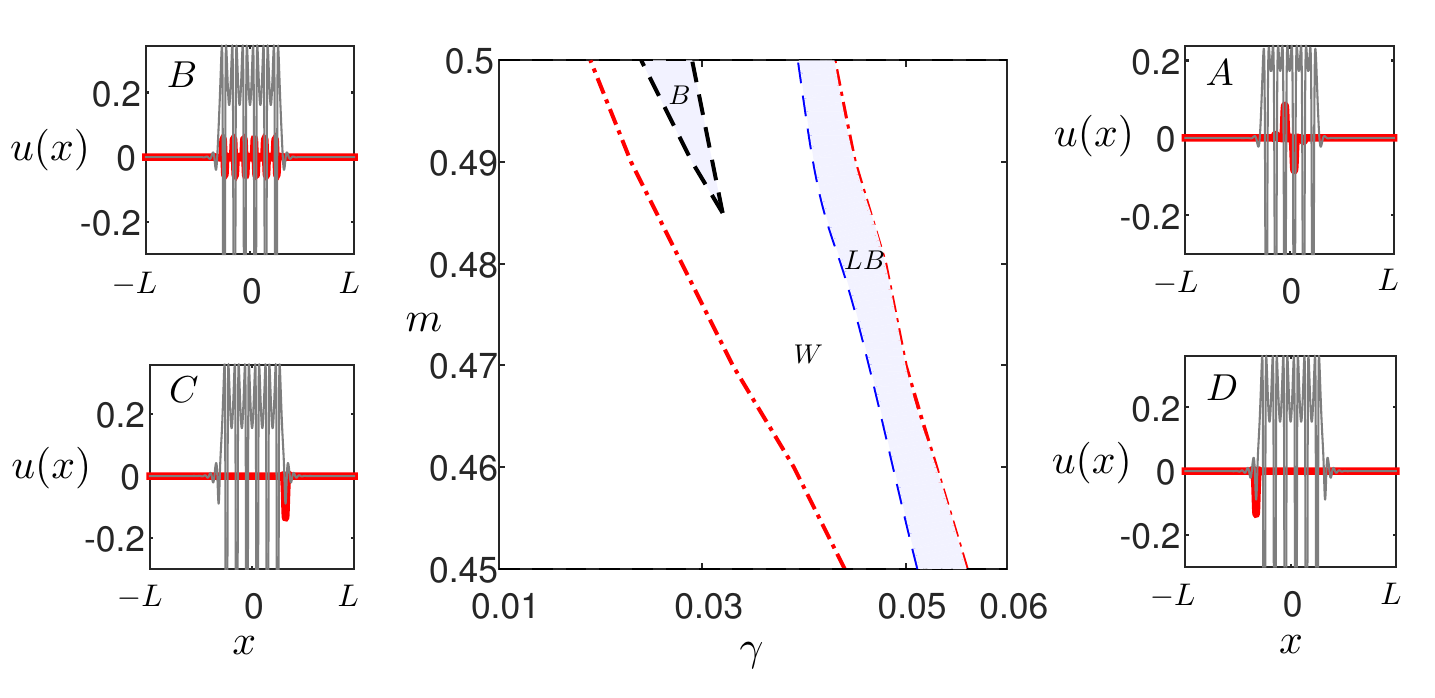}
	\caption{Regions of distinct transverse instabilities for localized stripe patterns according to~\eqref{eq:eignevalue_problem_2D_stripes} and the respective eigenfunctions: Localized--body ($LB$), body ($B$), wall ($W$) and the profiles $A$, $B$, and $C,D$, receptively. The dash-dotted line marks the limits of the snaking region. The analysis corresponds the $L_0$ branch between 5 and 6 saddle nodes. Profiles ($A$--$D$) show the most unstable eigenfunction in the background of the 1D localized state profile, as obtained numerically on a domain $2L=128$ that is much larger than the typical critical period $2 \pi/k_c \sim O(1)$. Parameters: $a=\tau=0$ and $\sigma=1.$}\label{fig:FC}
\end{figure}

\begin{figure}[tp]
\includegraphics[width=\textwidth]{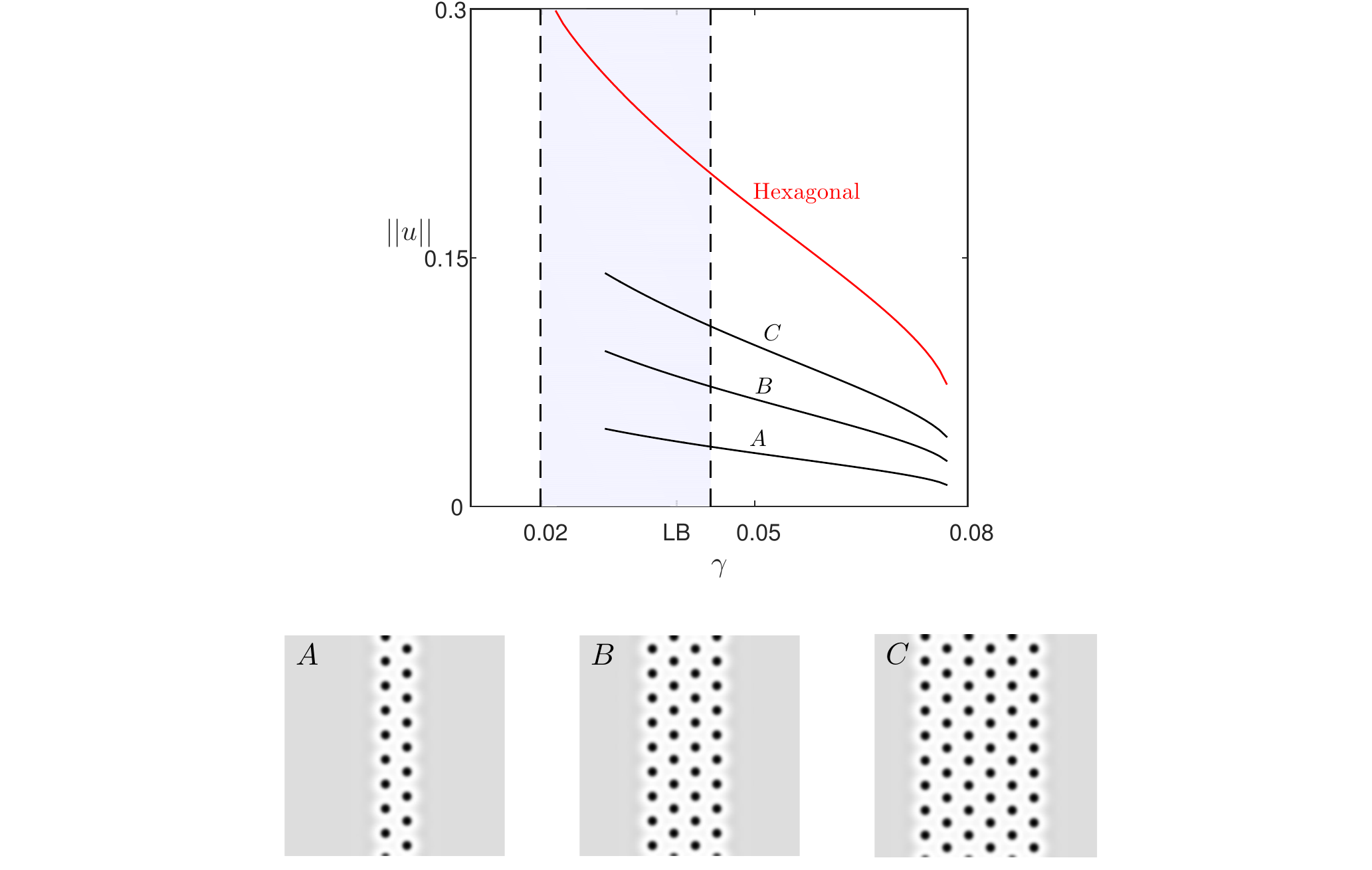}
	\caption{Top panel: Stable branches of localized hexagons and a spatially extended hexagonal pattern (light line); LB marks the left boundary of the localized body modes (see Fig.~\ref{fig:FC}). Bottom panel: Typical profiles of localized solutions along the branches that are respectively marked in the top panel.
	The localized hexagons have been computed on [$40\times 40$] domains while spatially extended hexagons on [$5\frac{2\pi}{k_c}\times 5\frac{2\pi}{k_c}$] domains, where $k_c\simeq 2.828$. Parameters: $m=0.5$, $a=\tau=0$, and $\sigma=1$.}\label{fig:branch_hex}
\end{figure}

The growth rates $\beta_{max}(k)$, show distinct instability mechanisms~\cite{burke2007homoclinic}, which  are consistent with the emergence of hexagonal/spot type solutions~\cite{choksi20112d}:
\begin{itemize}
	\item  The primary instability is of long wavelength (zigzag) type, as shown in Figure~\ref{fig:dispersion_wall_body}(a) and respectively by direct numerical integration in Figure~\ref{fig:TI}(a). Consequently, we refer to this mode as a \textit{body} mode;
	
	\item As $\gamma$ is varied, there is a growth of a secondary finite mode, to which we refer as a \textit{wall} mode, as shown in Figure~\ref{fig:dispersion_wall_body}(b) and Figure~\ref{fig:TI}(b).
	
	\item In addition, there exists a \textit{localized--body} mode that corresponds to spatially local domain within the interior of stripes, as shown in Figure~\ref{fig:TI}(c).
	
\end{itemize}
The perturbation that is associated with the body instability develops on a slower time scale than the wall and localized--body instabilities, as also indicated by larger growth rate ($\beta_{max}$) values in Figure~\ref{fig:dispersion_wall_body}. A summary of the transverse instabilities and the corresponding eigenfunctions is shown in the parameter space ($\gamma,m$), see Figure~\ref{fig:FC}. Although not in the scope of this paper, localized spots can also persist in this type of problems~\cite{lloyd2008localized,choksi20112d} due to equal energy between spots and a uniform state. In Figure \ref{fig:branch_hex}, we present several coexisting stable branches of localized hexagonal patterns which were computed by direct numerical integration of~\eqref{eq:intro_pde1}; the shaded region denotes 1D homoclinic snaking region.


\begin{figure}[tp]
	\large{ (a)\includegraphics[width=0.45\textwidth]{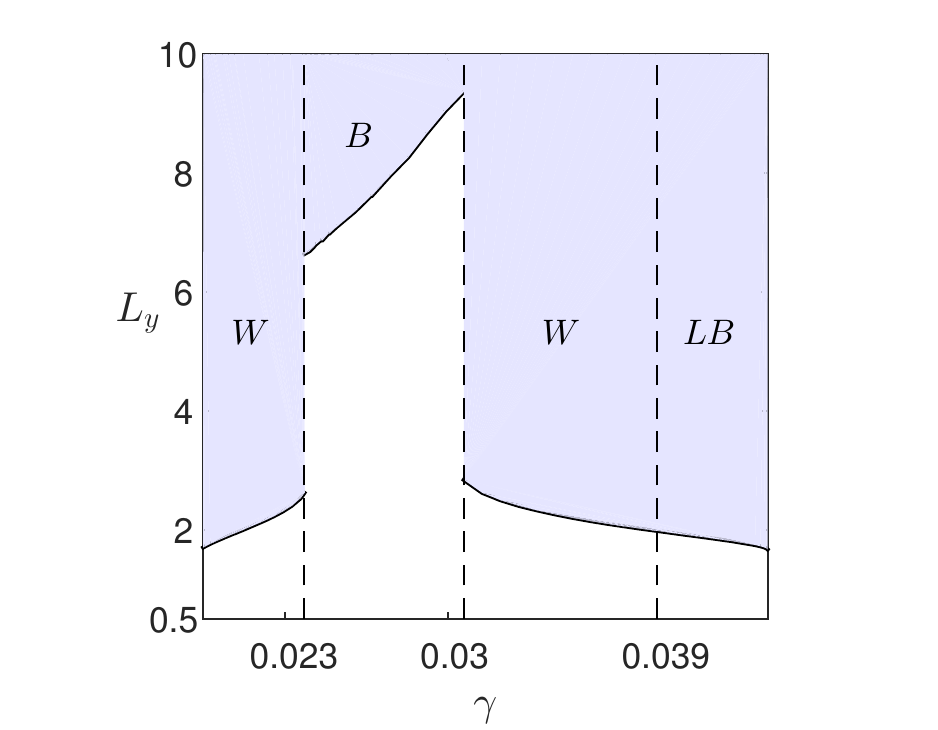}
	(b)}\includegraphics[width=0.45\textwidth]{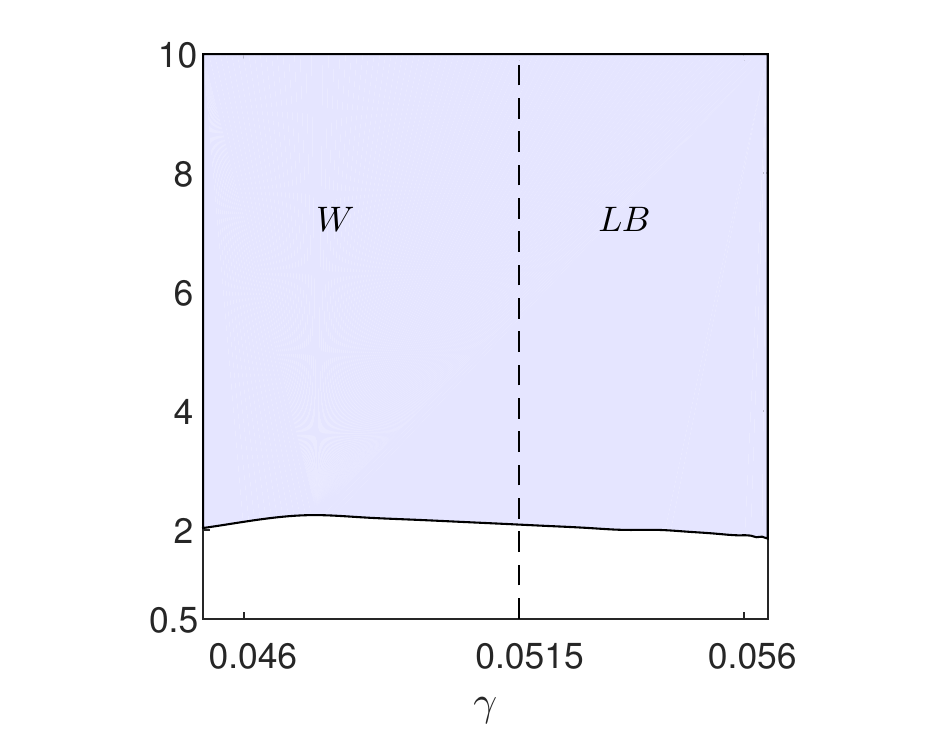}
	\caption{Dependence of transverse instabilities on domain size in $y$ direction within the snaking region, i.e., numerical results obtained via~\eqref{eq:eignevalue_problem_2D_stripes} showing the critical domain size according to the instability type. The computations have been performed with respect to cuts in Figure~\ref{fig:FC} at (a) $m=0.5$, (b) $m=0.45$. Other parameters: $a=\tau=0$, $\sigma=1$.}\label{fig:FC_Ly}
\end{figure}

Additionally to dependence on equation parameters, spatially localized stripe instability depends on a finite domain size in $y$ direction, i.e., on $L_y$. The latter is attributed to the preferable energy between the hexagonal and the stripe patterns. Thus, in the region where body modes exist (larger $m$ that corresponds to asymmetry), the stability of localized stripes over extended $L_y$ domains can be achieved due to the long wavenumber instability mechanism, as shown by the shaded region in Figure~\ref{fig:FC_Ly}(a). Otherwise in the sole presence of localized eigenfunctions (i.e., for wall and localized--body), localized stripes are obtained for relatively small domains ($L_y<\lambda_{H}$), where $\lambda_{H}$ the hexagonal characteristic period, as shown in Figure~\ref{fig:FC_Ly}(b).

\subsection{Extended Ohta–Kawasaki model, $a,\tau\neq 0$}	\label{sec:Loc_EOK}

\begin{figure}[tp]
\includegraphics[width=\textwidth]{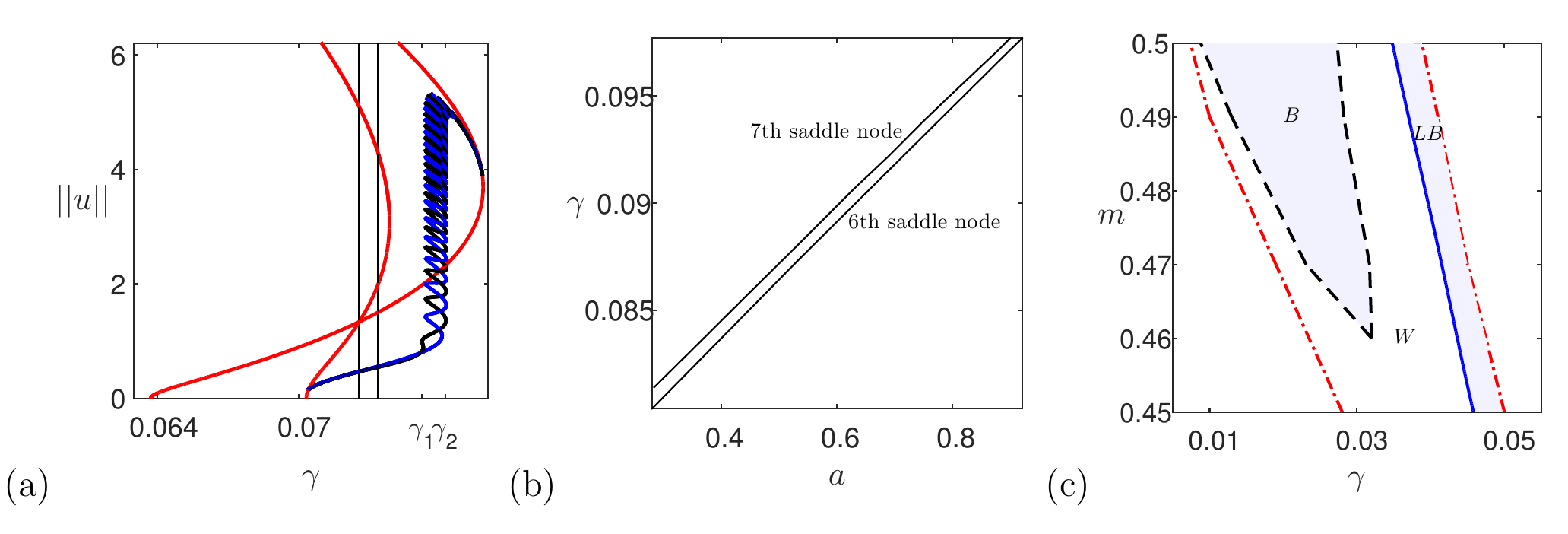}
	\caption{(a) Homoclinic snaking region ($\gamma_1<\gamma<\gamma_2$) for the EOK equation ($a,\tau\neq 0$) obtained via numerical continuation, see Figure~\ref{fig:snaking_ok}, where $\gamma_1\simeq0.0752$, $\gamma_2\simeq 0.0762$. Two branches of localized states are indicated by dark intertwined lines while light lines indicate periodic states. Stability is the same as in Figure~\ref{fig:snaking_ok} and not shown here. The localized states terminate on a periodic branch $\frac{35}{40}k_c$, where $k_c\simeq 1.923$. For comparison, we indicate the snaking region at $a=0$ (vertical light lines).  Parameters: $a=0.1$, $m=0.4$, $\tau=0$, $\sigma=1$. (b) Persistence of the snaking region in parameter space ($a,\gamma$). The snaking region lies in between the two lines that correspond to the sixth and seventh saddle nodes in (a). The locations have been obtained numerically via continuation of saddle nodes. (c) Regions of distinct transverse instabilities as in Figure~\ref{fig:FC} but for $\tau=-0.1$, $a=0.1$.
	}\label{fig:snaking_eps_non_tau}
\end{figure}

According to the linear and weakly nonlinear analyses, $a\ne 0$ and $\tau\ne0$ do not indicate any qualitative differences regarding the presence of localized states, see Section~\ref{sec:anal}. Numerical continuation shows that both 1D homoclinic snaking and the respective selection mechanism, are indeed present although being shifted, as shown in Figure~\ref{fig:snaking_eps_non_tau}(a). To generalize the results for a larger range of parameters, we chose a pair of adjacent saddle nodes (specifically sixth and seventh) and obtained their locations and consequently the existence of homoclinic snaking region, see Figure~\ref{fig:snaking_eps_non_tau}(b). The 2D stability properties of localized stripes are also persist, although as expected with some shift in parameter range, as shown in Figure~\ref{fig:snaking_eps_non_tau}(c). Specifically, it is shown that asymmetry in charged domains increases the region of stability for localized stripes. To this end, the localized patterns and mechanisms obtained for the classical OK model ($a=\tau=0$) hold also for the more general case and due to simplicity can be exploited for future studies.

\section{Discussion} \label{sec:concl}	

Self--assembly of spatially localized states in one- and two-space dimensions has been studied in a context of electrically charged and mass conserved media (i.e., a system which exhibits both short and long range interactions). As a case study system, we focus on the Ohta-Kawasaki framework which incorporates, in addition, difference in the properties of the charged domains and non-uniform permittivity. While the OK model has been studied extensively in the context of spatially extended patterns that arise through super--critical bifurcations, it is the sub--critical bifurcations that give rise to a rich behavior of localized states. Although the OK system is conservative, the self--assembly behaves differently than typical phase field model (a.k.a. conserved Swift-Hohenberg)~\cite{beaume2013convectons,dawes2008localized,thiele2013localized}. In particular, the OK system obeys the linear stability properties (dispersion relation) of dissipative systems, such as Swift-Hohenberg~\cite{burke2007homoclinic} and Gierer-Meinhardt~\cite{yochelis2008formation} and the formed localized states in 1D create a vertical (in contrast to  slanted, as in~\cite{beaume2013convectons,dawes2008localized,thiele2013localized}) homoclinic snaking structure that corresponds to heteroclinic cycles linking uniform and periodic solutions. 

The reason for the distinct behavior is attributed to the coupling with non--local interactions that precludes unique chemical potential for both the mass and the electrical field. Without the electrical field, any change in the model control parameter corresponds to a respective change in chemical potential~\cite{thiele2013localized}, which leads in turn to slanted snaking. Although the EOK model also posses the property of energy minimization under mass constraint, the first variation of the energy depends also on the solution of Poisson's equation~\eqref{eq:poisson}. However, the latter is not unique under Neumann boundary conditions and sets a degree of freedom in setting one of the Lagrange multipliers, see~\eqref{eq:Lag} and the appendix.  Consequenly, vertical snaking can form without violation of thermodynamic principles. 

The impact of Coulombic interactions is reflected also in the asymptotic analysis and the form of the amplitude equations: Poissons' equation~\eqref{eq:poisson} is solvable only when $\int_{\Omega} u-m \ d\Omega=0$.  Therefore mass of EOK~\eqref{eq:intro_pde1} solutions cannot be arbitrarily determined via the initial condition, as occurs, for example in cSH~\cite{thiele2013localized}. The arbitrary choice of mass in cSH, is reflected in the dispersion relation via the ever neutrality of $\lambda(k=0)=0$~\cite{golovin2003self,dawes2008localized}, different from the EOK model~\eqref{eq:intro_pde1} where $\lambda(k=0)<0$, see Figure \ref{fig:disper}. Coupling to $k=0$ results in two coupled amplitude equations with nonlocal properties~\cite{dawes2008localized}.
Consequently, our study gives rise to a counter-intuitive picture that while conservation in local equations leads to non-local effects in the amplitude equations, in our case coupling to non-locality leads to local effects in the amplitude equations, bearing similarity to dissipative systems~\cite{knobloch2015spatial}.

Localized stripes which are the natural extension in 2D have been also investigated. We found distinct instability types associated with spatially extended eigenfunctions (body modes) and localized eigenfunctions (wall and body--localized modes)~\cite{burke2007homoclinic}, as demonstrated in Figures~\ref{fig:TI} and~\ref{fig:FC}. Specifically, the instability via localized eigenfunctions leads to localized spots while body eigenfunctions lead to stripes embedded in a hexagonal background; stability of stripes to body modes persists over a larger domain as compared with the localized modes (see Figure~\ref{fig:FC_Ly}). In a more general context, localized pattern formation in the vicinity of Maxwell points, e.g., an energy equivalence between periodic and  uniform states, allows to interpret existing phase maps in 2D or 3D to map possible localized structures in high-dimensions and thus applications. For example, the phase diagram devised for the OK equation~\cite{choksi20112d}, discloses an interface between stable periodic spots and strips. At the interface, the two states should presumably have the same energy and therefore, it is likely that localized strips on a spotted background would also coexist. On the other hand, in the absence of interface between stable stripes and uniform (disordered) states, it is unlikely that localized strips would be present. Indeed, both observations are consistent with our results. To this end, we believe that insights into the localized self--assembly mechanism developed here for conserved and charge systems, can be exploited to design by simple means, low cost isolated morphologies that can be attractive for example, to nano electronics~\cite{park2003enabling,stoykovich2007directed,kim2009block,farrell2010self}.

\section*{Acknowledgment}
This research was done in the framework of the Grand Technion Energy Program (GTEP) and of the BGU Energy Initiative Program, and supported by the Adelis Foundation for renewable energy research. N.G. also acknowledges the support from the Technion VPR fund and from EU Marie--Curie CIG Grant 2018620.

	\appendix*
	\section{Derivation of Newell-Whitehead-Segel amplitude equation} \label{sec:appendixa}
	In the vicinity of $\gamma_c$, we seek asymptotic solutions for~$\epsilon\ll1$ with multiple scales of the form
	\begin{subequations} \label{eq:expansion}
		\begin{equation}\label{eq:expansiona}
		u(t,x,y)=\sqrt{\epsilon}u_1(t,T,x,X,y,Y)+\epsilon u_2(t,T,x,X,y,Y)+\cdots,
		\end{equation}
		\begin{equation}\label{eq:expansionb}
		\phi(x,y)=\sqrt{\epsilon}\phi_1(x,X,y,Y)+\epsilon \phi_2(x,X,y,Y)+\cdots.
		\end{equation}
		Substituting~\eqref{eq:expansion} in~\eqref{eq:intro_pde1}, and matching terms in each order we obtain: 
	\end{subequations}\\\\
	$O(\sqrt{\epsilon}):$
	\begin{equation}\label{eq:sqrt epsilon order for u}
	\begin{split}
	0=&L(u_1,\phi_1):=-\gamma_c\partial_{x}^4u_1-(1-3m^2+2\tau m)\partial_{x}^2u_1+\sigma \partial_{x}^2\phi_1,\\
	0=&G(u_1,\phi_1):=(am+1)\partial_{x}^2\phi_1+u_1,
	\end{split}
	\end{equation}
	$O(\epsilon):$
	\begin{equation}\label{eq: epsilon order for u}
	\begin{split}
	0=&L(u_2,\phi_2)-4\gamma_c \partial_{xxxX}u_1-2(1-3m^2+2\tau m)\partial_{xX}u_1-2\gamma_c \partial_{xxYY}u_1 \\
	&+(3m-\tau)\partial_{xx}(u_1u_2)-(1-3m^2+2\tau m)\partial_{YY}u_1+2\sigma\partial_{xX}\phi_1+\sigma \partial_{YY}\phi_1,\\
	0=&G(u_2,\phi_2)+a\partial_x u_1\cdot \partial_x \phi_1+au_1\cdot\partial_x^2\phi_1+
	2(am+1)\partial_{xX}\phi_1+(am+1)\partial_{YY}\phi_1.
	\end{split}
	\end{equation}
	$O(\epsilon^{\frac{3}{2}}):$
	\begin{equation}\label{eq:epsilon1.5a}
	\begin{split}
	\partial_T u_1=&L(u_3,\phi_3)+\mathcal{N}_u[u_1,u_2,\phi_1,\phi_2],\\ 0=&G(u_3,\phi_3)+\mathcal{N}_\phi[u_1,u_2,\phi_1,\phi_2],
	\end{split}
	\end{equation}
	where~$\mathcal{N}_u$ and~$\mathcal{N}_\phi$ are nonlinear operators, which for brevity are not explicitly presented.
	Neumann boundary conditions are imposed for~$(u_i,\phi_i)$ at all orders.
	The solution of the first-order equation~\eqref{eq:sqrt epsilon order for u} is
	\begin{equation}\label{eq:sol_sqrt_epsilon}
	\begin{split}
	u_1=&A(X,Y,T)e^{ik_c x}+c.c.,\qquad \phi_1=Z(X,Y)+\frac{A(X,Y,T)}{k_c^2(am+1)}e^{ik_c x}+c.c.
	\end{split}
	\end{equation}
	The solution of the second-order equation~\eqref{eq: epsilon order for u} is
	\begin{equation}~\label{eq:sol_epsilon}
	\begin{split}
	u_2=&A_3(X,Y,T)e^{ik_c x}+\frac{2}{9}\cdot\frac{15am^2-6am\tau -a+12m-4\tau}
	{3am^3-2am^2\tau-am+3m^2-2m\tau-1}A^2(X,Y,T)e^{2ik_c x}+c.c.,\\
	\phi_2=&F(X,Y)+M(X,Y,T)e^{ik_c x}+\frac{1}{9}\cdot\frac{3am^2-3am\tau-2a-3m+\tau}
	{\sigma(am+1)}A^2(X,Y,T)e^{2ik_c x}+c.c.
	\end{split}
	\end{equation}
	Substituting the solutions~\eqref{eq:sol_sqrt_epsilon} and~\eqref{eq:sol_epsilon} into the third-order equation~\eqref{eq:epsilon1.5a} yields the form
	\begin{equation}\label{eq:third_order_eq}
	\begin{split}
	A_T e^{ik_c x}&=\mathcal{N}_0[u,\phi]+\mathcal{N}_1[u,\phi]e^{ik_cx}+\mathcal{N}_2[u,\phi] e^{2ik_cx},
	\end{split}
	\end{equation}
	where~$\mathcal{N}_i$ are operators of~$[u_1,u_2,u_3,\phi_1,\phi_2,\phi_3]$. Orthogonality of the Fourier modes implies that the coefficients of~$e^{ik_c x}$ in~\eqref{eq:third_order_eq} equate, yielding
	\begin{equation}
	A_T=\mathcal{N}_1[A]=\hat{\gamma}k_c^4A+f|A|^2A+\left(\sqrt{2\eta}\frac{\partial}{\partial X}-i\sqrt{\gamma_c}\frac{\partial^2}{\partial Y^2} \right)^2A.
	\end{equation}
	Similarly, the equations for $Z(X,Y)$ in~\eqref{eq:sol_sqrt_epsilon} read as
	$$
	Z_{XX}=Z_{YY}=0,\qquad \left. \frac{\partial Z}{\partial X}\right|_{\partial \Omega}=\left. \frac{\partial Z}{\partial Y}\right|_{\partial \Omega}=0.
	$$
	Therefore, $Z$ is a constant function.  This implies that, as expected,~$\phi$ is not unique under Neumann boundary conditions, see~\eqref{eq:sol_sqrt_epsilon}, but rather determined up to a constant.


\begin{thebibliography}{36}%
\makeatletter
\providecommand \@ifxundefined [1]{%
 \@ifx{#1\undefined}
}%
\providecommand \@ifnum [1]{%
 \ifnum #1\expandafter \@firstoftwo
 \else \expandafter \@secondoftwo
 \fi
}%
\providecommand \@ifx [1]{%
 \ifx #1\expandafter \@firstoftwo
 \else \expandafter \@secondoftwo
 \fi
}%
\providecommand \natexlab [1]{#1}%
\providecommand \enquote  [1]{``#1''}%
\providecommand \bibnamefont  [1]{#1}%
\providecommand \bibfnamefont [1]{#1}%
\providecommand \citenamefont [1]{#1}%
\providecommand \href@noop [0]{\@secondoftwo}%
\providecommand \href [0]{\begingroup \@sanitize@url \@href}%
\providecommand \@href[1]{\@@startlink{#1}\@@href}%
\providecommand \@@href[1]{\endgroup#1\@@endlink}%
\providecommand \@sanitize@url [0]{\catcode `\\12\catcode `\$12\catcode
  `\&12\catcode `\#12\catcode `\^12\catcode `\_12\catcode `\%12\relax}%
\providecommand \@@startlink[1]{}%
\providecommand \@@endlink[0]{}%
\providecommand \url  [0]{\begingroup\@sanitize@url \@url }%
\providecommand \@url [1]{\endgroup\@href {#1}{\urlprefix }}%
\providecommand \urlprefix  [0]{URL }%
\providecommand \Eprint [0]{\href }%
\providecommand \doibase [0]{http://dx.doi.org/}%
\providecommand \selectlanguage [0]{\@gobble}%
\providecommand \bibinfo  [0]{\@secondoftwo}%
\providecommand \bibfield  [0]{\@secondoftwo}%
\providecommand \translation [1]{[#1]}%
\providecommand \BibitemOpen [0]{}%
\providecommand \bibitemStop [0]{}%
\providecommand \bibitemNoStop [0]{.\EOS\space}%
\providecommand \EOS [0]{\spacefactor3000\relax}%
\providecommand \BibitemShut  [1]{\csname bibitem#1\endcsname}%
\let\auto@bib@innerbib\@empty
\bibitem [{\citenamefont {Chidichimo}\ and\ \citenamefont
  {Filippelli}(2010)}]{ChFi:2010}%
  \BibitemOpen
  \bibfield  {author} {\bibinfo {author} {\bibfnamefont {G.}~\bibnamefont
  {Chidichimo}}\ and\ \bibinfo {author} {\bibfnamefont {L.}~\bibnamefont
  {Filippelli}},\ }\href@noop {} {\bibfield  {journal} {\bibinfo  {journal}
  {International Journal of Photoenergy}\ ,\ \bibinfo {pages} {123534}}
  (\bibinfo {year} {2010})}\BibitemShut {NoStop}%
\bibitem [{\citenamefont {Wolden}\ \emph {et~al.}(2011)\citenamefont {Wolden},
  \citenamefont {Kurtin}, \citenamefont {Baxter}, \citenamefont {Repins},
  \citenamefont {Shaheen}, \citenamefont {Torvik}, \citenamefont {Rockett},
  \citenamefont {Fthenakis},\ and\ \citenamefont {Aydil}}]{OPV_Rev:2011}%
  \BibitemOpen
  \bibfield  {author} {\bibinfo {author} {\bibfnamefont {C.~A.}\ \bibnamefont
  {Wolden}}, \bibinfo {author} {\bibfnamefont {J.}~\bibnamefont {Kurtin}},
  \bibinfo {author} {\bibfnamefont {J.~B.}\ \bibnamefont {Baxter}}, \bibinfo
  {author} {\bibfnamefont {I.}~\bibnamefont {Repins}}, \bibinfo {author}
  {\bibfnamefont {S.~E.}\ \bibnamefont {Shaheen}}, \bibinfo {author}
  {\bibfnamefont {J.~T.}\ \bibnamefont {Torvik}}, \bibinfo {author}
  {\bibfnamefont {A.~A.}\ \bibnamefont {Rockett}}, \bibinfo {author}
  {\bibfnamefont {V.~M.}\ \bibnamefont {Fthenakis}}, \ and\ \bibinfo {author}
  {\bibfnamefont {E.~S.}\ \bibnamefont {Aydil}},\ }\href@noop {} {\bibfield
  {journal} {\bibinfo  {journal} {Journal of Vacuum Science \& Technology A}\
  }\textbf {\bibinfo {volume} {29}},\ \bibinfo {pages} {030801} (\bibinfo
  {year} {2011})}\BibitemShut {NoStop}%
\bibitem [{\citenamefont {Promislow}\ and\ \citenamefont
  {Wetton}(2009)}]{promislow2009pem}%
  \BibitemOpen
  \bibfield  {author} {\bibinfo {author} {\bibfnamefont {K.}~\bibnamefont
  {Promislow}}\ and\ \bibinfo {author} {\bibfnamefont {B.}~\bibnamefont
  {Wetton}},\ }\href@noop {} {\bibfield  {journal} {\bibinfo  {journal} {SIAM
  Journal on Applied Mathematics}\ }\textbf {\bibinfo {volume} {70}},\ \bibinfo
  {pages} {369} (\bibinfo {year} {2009})}\BibitemShut {NoStop}%
\bibitem [{\citenamefont {Mauritz}\ and\ \citenamefont
  {Moore}(2004)}]{mauritz2004state}%
  \BibitemOpen
  \bibfield  {author} {\bibinfo {author} {\bibfnamefont {K.~A.}\ \bibnamefont
  {Mauritz}}\ and\ \bibinfo {author} {\bibfnamefont {R.~B.}\ \bibnamefont
  {Moore}},\ }\href@noop {} {\bibfield  {journal} {\bibinfo  {journal}
  {Chemical Reviews}\ }\textbf {\bibinfo {volume} {104}},\ \bibinfo {pages}
  {4535} (\bibinfo {year} {2004})}\BibitemShut {NoStop}%
\bibitem [{\citenamefont {Tsori}(2009)}]{tsori2009polymers}%
  \BibitemOpen
  \bibfield  {author} {\bibinfo {author} {\bibfnamefont {Y.}~\bibnamefont
  {Tsori}},\ }\href@noop {} {\emph {\bibinfo {title} {Polymers, liquids and
  colloids in electric fields: interfacial instabilities, orientation and phase
  transitions}}},\ Vol.~\bibinfo {volume} {2}\ (\bibinfo  {publisher} {World
  Scientific},\ \bibinfo {year} {2009})\BibitemShut {NoStop}%
\bibitem [{\citenamefont {Chen}(2002)}]{chen2002phase}%
  \BibitemOpen
  \bibfield  {author} {\bibinfo {author} {\bibfnamefont {L.-Q.}\ \bibnamefont
  {Chen}},\ }\href@noop {} {\bibfield  {journal} {\bibinfo  {journal} {Annual
  review of materials research}\ }\textbf {\bibinfo {volume} {32}},\ \bibinfo
  {pages} {113} (\bibinfo {year} {2002})}\BibitemShut {NoStop}%
\bibitem [{\citenamefont {Golovin}\ \emph {et~al.}(2003)\citenamefont
  {Golovin}, \citenamefont {Davis},\ and\ \citenamefont
  {Voorhees}}]{golovin2003self}%
  \BibitemOpen
  \bibfield  {author} {\bibinfo {author} {\bibfnamefont {A.}~\bibnamefont
  {Golovin}}, \bibinfo {author} {\bibfnamefont {S.}~\bibnamefont {Davis}}, \
  and\ \bibinfo {author} {\bibfnamefont {P.}~\bibnamefont {Voorhees}},\
  }\href@noop {} {\bibfield  {journal} {\bibinfo  {journal} {Physical Review
  E}\ }\textbf {\bibinfo {volume} {68}},\ \bibinfo {pages} {056203} (\bibinfo
  {year} {2003})}\BibitemShut {NoStop}%
\bibitem [{\citenamefont {Yochelis}\ \emph {et~al.}(2015)\citenamefont
  {Yochelis}, \citenamefont {Singh},\ and\ \citenamefont
  {Visoly-Fisher}}]{yochelis2015coupling}%
  \BibitemOpen
  \bibfield  {author} {\bibinfo {author} {\bibfnamefont {A.}~\bibnamefont
  {Yochelis}}, \bibinfo {author} {\bibfnamefont {M.~B.}\ \bibnamefont {Singh}},
  \ and\ \bibinfo {author} {\bibfnamefont {I.}~\bibnamefont {Visoly-Fisher}},\
  }\href@noop {} {\bibfield  {journal} {\bibinfo  {journal} {Chemistry of
  Materials}\ }\textbf {\bibinfo {volume} {27}},\ \bibinfo {pages} {4169}
  (\bibinfo {year} {2015})}\BibitemShut {NoStop}%
\bibitem [{\citenamefont {Ohta}\ and\ \citenamefont
  {Kawasaki}(1986)}]{ohta1986equilibrium}%
  \BibitemOpen
  \bibfield  {author} {\bibinfo {author} {\bibfnamefont {T.}~\bibnamefont
  {Ohta}}\ and\ \bibinfo {author} {\bibfnamefont {K.}~\bibnamefont
  {Kawasaki}},\ }\href@noop {} {\bibfield  {journal} {\bibinfo  {journal}
  {Macromolecules}\ }\textbf {\bibinfo {volume} {19}},\ \bibinfo {pages} {2621}
  (\bibinfo {year} {1986})}\BibitemShut {NoStop}%
\bibitem [{\citenamefont {Choksi}\ \emph {et~al.}(2009)\citenamefont {Choksi},
  \citenamefont {Peletier},\ and\ \citenamefont {Williams}}]{choksi2009phase}%
  \BibitemOpen
  \bibfield  {author} {\bibinfo {author} {\bibfnamefont {R.}~\bibnamefont
  {Choksi}}, \bibinfo {author} {\bibfnamefont {M.~A.}\ \bibnamefont
  {Peletier}}, \ and\ \bibinfo {author} {\bibfnamefont {J.}~\bibnamefont
  {Williams}},\ }\href@noop {} {\bibfield  {journal} {\bibinfo  {journal} {SIAM
  Journal on Applied Mathematics}\ }\textbf {\bibinfo {volume} {69}},\ \bibinfo
  {pages} {1712} (\bibinfo {year} {2009})}\BibitemShut {NoStop}%
\bibitem [{\citenamefont {Gavish}\ and\ \citenamefont
  {Yochelis}(2016)}]{gavish2016theory}%
  \BibitemOpen
  \bibfield  {author} {\bibinfo {author} {\bibfnamefont {N.}~\bibnamefont
  {Gavish}}\ and\ \bibinfo {author} {\bibfnamefont {A.}~\bibnamefont
  {Yochelis}},\ }\href@noop {} {\bibfield  {journal} {\bibinfo  {journal}
  {Journal of Physical Chemistry Letters}\ }\textbf {\bibinfo {volume} {7}},\
  \bibinfo {pages} {1121} (\bibinfo {year} {2016})}\BibitemShut {NoStop}%
\bibitem [{\citenamefont {Knobloch}(2015)}]{knobloch2015spatial}%
  \BibitemOpen
  \bibfield  {author} {\bibinfo {author} {\bibfnamefont {E.}~\bibnamefont
  {Knobloch}},\ }\href@noop {} {\bibfield  {journal} {\bibinfo  {journal}
  {conmatphys}\ }\textbf {\bibinfo {volume} {6}},\ \bibinfo {pages} {325}
  (\bibinfo {year} {2015})}\BibitemShut {NoStop}%
\bibitem [{\citenamefont {Choksi}\ \emph {et~al.}(2011)\citenamefont {Choksi},
  \citenamefont {Maras},\ and\ \citenamefont {Williams}}]{choksi20112d}%
  \BibitemOpen
  \bibfield  {author} {\bibinfo {author} {\bibfnamefont {R.}~\bibnamefont
  {Choksi}}, \bibinfo {author} {\bibfnamefont {M.}~\bibnamefont {Maras}}, \
  and\ \bibinfo {author} {\bibfnamefont {J.}~\bibnamefont {Williams}},\
  }\href@noop {} {\bibfield  {journal} {\bibinfo  {journal} {SIAM Journal on
  Applied Dynamical Systems}\ }\textbf {\bibinfo {volume} {10}},\ \bibinfo
  {pages} {1344} (\bibinfo {year} {2011})}\BibitemShut {NoStop}%
\bibitem [{\citenamefont {Choksi}\ and\ \citenamefont
  {Sternberg}(2006)}]{choksi2006periodic}%
  \BibitemOpen
  \bibfield  {author} {\bibinfo {author} {\bibfnamefont {R.}~\bibnamefont
  {Choksi}}\ and\ \bibinfo {author} {\bibfnamefont {P.}~\bibnamefont
  {Sternberg}},\ }\href@noop {} {\bibfield  {journal} {\bibinfo  {journal}
  {Interfaces and Free Boundaries}\ }\textbf {\bibinfo {volume} {8}},\ \bibinfo
  {pages} {371} (\bibinfo {year} {2006})}\BibitemShut {NoStop}%
\bibitem [{\citenamefont {Park}\ \emph {et~al.}(2003)\citenamefont {Park},
  \citenamefont {Yoon},\ and\ \citenamefont {Thomas}}]{park2003enabling}%
  \BibitemOpen
  \bibfield  {author} {\bibinfo {author} {\bibfnamefont {C.}~\bibnamefont
  {Park}}, \bibinfo {author} {\bibfnamefont {J.}~\bibnamefont {Yoon}}, \ and\
  \bibinfo {author} {\bibfnamefont {E.~L.}\ \bibnamefont {Thomas}},\
  }\href@noop {} {\bibfield  {journal} {\bibinfo  {journal} {Polymer}\ }\textbf
  {\bibinfo {volume} {44}},\ \bibinfo {pages} {6725} (\bibinfo {year}
  {2003})}\BibitemShut {NoStop}%
\bibitem [{\citenamefont {Stoykovich}\ \emph {et~al.}(2007)\citenamefont
  {Stoykovich}, \citenamefont {Kang}, \citenamefont {Daoulas}, \citenamefont
  {Liu}, \citenamefont {Liu}, \citenamefont {de~Pablo}, \citenamefont
  {M{\"u}ller},\ and\ \citenamefont {Nealey}}]{stoykovich2007directed}%
  \BibitemOpen
  \bibfield  {author} {\bibinfo {author} {\bibfnamefont {M.~P.}\ \bibnamefont
  {Stoykovich}}, \bibinfo {author} {\bibfnamefont {H.}~\bibnamefont {Kang}},
  \bibinfo {author} {\bibfnamefont {K.~C.}\ \bibnamefont {Daoulas}}, \bibinfo
  {author} {\bibfnamefont {G.}~\bibnamefont {Liu}}, \bibinfo {author}
  {\bibfnamefont {C.-C.}\ \bibnamefont {Liu}}, \bibinfo {author} {\bibfnamefont
  {J.~J.}\ \bibnamefont {de~Pablo}}, \bibinfo {author} {\bibfnamefont
  {M.}~\bibnamefont {M{\"u}ller}}, \ and\ \bibinfo {author} {\bibfnamefont
  {P.~F.}\ \bibnamefont {Nealey}},\ }\href@noop {} {\bibfield  {journal}
  {\bibinfo  {journal} {Acs Nano}\ }\textbf {\bibinfo {volume} {1}},\ \bibinfo
  {pages} {168} (\bibinfo {year} {2007})}\BibitemShut {NoStop}%
\bibitem [{\citenamefont {Kim}\ \emph {et~al.}(2009)\citenamefont {Kim},
  \citenamefont {Park},\ and\ \citenamefont {Hinsberg}}]{kim2009block}%
  \BibitemOpen
  \bibfield  {author} {\bibinfo {author} {\bibfnamefont {H.-C.}\ \bibnamefont
  {Kim}}, \bibinfo {author} {\bibfnamefont {S.-M.}\ \bibnamefont {Park}}, \
  and\ \bibinfo {author} {\bibfnamefont {W.~D.}\ \bibnamefont {Hinsberg}},\
  }\href@noop {} {\bibfield  {journal} {\bibinfo  {journal} {Chemical reviews}\
  }\textbf {\bibinfo {volume} {110}},\ \bibinfo {pages} {146} (\bibinfo {year}
  {2009})}\BibitemShut {NoStop}%
\bibitem [{\citenamefont {Farrell}\ \emph {et~al.}(2010)\citenamefont
  {Farrell}, \citenamefont {Petkov}, \citenamefont {Morris},\ and\
  \citenamefont {Holmes}}]{farrell2010self}%
  \BibitemOpen
  \bibfield  {author} {\bibinfo {author} {\bibfnamefont {R.~A.}\ \bibnamefont
  {Farrell}}, \bibinfo {author} {\bibfnamefont {N.}~\bibnamefont {Petkov}},
  \bibinfo {author} {\bibfnamefont {M.~A.}\ \bibnamefont {Morris}}, \ and\
  \bibinfo {author} {\bibfnamefont {J.~D.}\ \bibnamefont {Holmes}},\
  }\href@noop {} {\bibfield  {journal} {\bibinfo  {journal} {Journal of colloid
  and interface science}\ }\textbf {\bibinfo {volume} {349}},\ \bibinfo {pages}
  {449} (\bibinfo {year} {2010})}\BibitemShut {NoStop}%
\bibitem [{\citenamefont {Glasner}(2010)}]{glasner2010spatially}%
  \BibitemOpen
  \bibfield  {author} {\bibinfo {author} {\bibfnamefont {K.~B.}\ \bibnamefont
  {Glasner}},\ }\href@noop {} {\bibfield  {journal} {\bibinfo  {journal} {SIAM
  Journal on Applied Mathematics}\ }\textbf {\bibinfo {volume} {70}},\ \bibinfo
  {pages} {2045} (\bibinfo {year} {2010})}\BibitemShut {NoStop}%
\bibitem [{\citenamefont {Dawes}(2008)}]{dawes2008localized}%
  \BibitemOpen
  \bibfield  {author} {\bibinfo {author} {\bibfnamefont {J.~H.}\ \bibnamefont
  {Dawes}},\ }\href@noop {} {\bibfield  {journal} {\bibinfo  {journal} {SIAM
  Journal on Applied Dynamical Systems}\ }\textbf {\bibinfo {volume} {7}},\
  \bibinfo {pages} {186} (\bibinfo {year} {2008})}\BibitemShut {NoStop}%
\bibitem [{\citenamefont {Thiele}\ \emph {et~al.}(2013)\citenamefont {Thiele},
  \citenamefont {Archer}, \citenamefont {Robbins}, \citenamefont {Gomez},\ and\
  \citenamefont {Knobloch}}]{thiele2013localized}%
  \BibitemOpen
  \bibfield  {author} {\bibinfo {author} {\bibfnamefont {U.}~\bibnamefont
  {Thiele}}, \bibinfo {author} {\bibfnamefont {A.~J.}\ \bibnamefont {Archer}},
  \bibinfo {author} {\bibfnamefont {M.~J.}\ \bibnamefont {Robbins}}, \bibinfo
  {author} {\bibfnamefont {H.}~\bibnamefont {Gomez}}, \ and\ \bibinfo {author}
  {\bibfnamefont {E.}~\bibnamefont {Knobloch}},\ }\href@noop {} {\bibfield
  {journal} {\bibinfo  {journal} {Physical Review E}\ }\textbf {\bibinfo
  {volume} {87}},\ \bibinfo {pages} {042915} (\bibinfo {year}
  {2013})}\BibitemShut {NoStop}%
\bibitem [{\citenamefont {Orizaga}\ and\ \citenamefont
  {Glasner}(2016)}]{orizaga2016instability}%
  \BibitemOpen
  \bibfield  {author} {\bibinfo {author} {\bibfnamefont {S.}~\bibnamefont
  {Orizaga}}\ and\ \bibinfo {author} {\bibfnamefont {K.}~\bibnamefont
  {Glasner}},\ }\href@noop {} {\bibfield  {journal} {\bibinfo  {journal}
  {Physical Review E}\ }\textbf {\bibinfo {volume} {93}},\ \bibinfo {pages}
  {052504} (\bibinfo {year} {2016})}\BibitemShut {NoStop}%
\bibitem [{\citenamefont {Shiwa}(1997)}]{shiwa1997amplitude}%
  \BibitemOpen
  \bibfield  {author} {\bibinfo {author} {\bibfnamefont {Y.}~\bibnamefont
  {Shiwa}},\ }\href@noop {} {\bibfield  {journal} {\bibinfo  {journal} {Physics
  Letters A}\ }\textbf {\bibinfo {volume} {228}},\ \bibinfo {pages} {279}
  (\bibinfo {year} {1997})}\BibitemShut {NoStop}%
\bibitem [{\citenamefont {Champneys}(1998)}]{champneys1998homoclinic}%
  \BibitemOpen
  \bibfield  {author} {\bibinfo {author} {\bibfnamefont {A.~R.}\ \bibnamefont
  {Champneys}},\ }\href@noop {} {\bibfield  {journal} {\bibinfo  {journal}
  {Physica D: Nonlinear Phenomena}\ }\textbf {\bibinfo {volume} {112}},\
  \bibinfo {pages} {158} (\bibinfo {year} {1998})}\BibitemShut {NoStop}%
\bibitem [{\citenamefont {Hoyle}(2006)}]{hoyle2006pattern}%
  \BibitemOpen
  \bibfield  {author} {\bibinfo {author} {\bibfnamefont {R.~B.}\ \bibnamefont
  {Hoyle}},\ }\href@noop {} {\emph {\bibinfo {title} {Pattern formation: an
  introduction to methods}}}\ (\bibinfo  {publisher} {Cambridge University
  Press},\ \bibinfo {year} {2006})\BibitemShut {NoStop}%
\bibitem [{\citenamefont {Chapman}\ and\ \citenamefont
  {Kozyreff}(2009)}]{chapman2009exponential}%
  \BibitemOpen
  \bibfield  {author} {\bibinfo {author} {\bibfnamefont {S.}~\bibnamefont
  {Chapman}}\ and\ \bibinfo {author} {\bibfnamefont {G.}~\bibnamefont
  {Kozyreff}},\ }\href@noop {} {\bibfield  {journal} {\bibinfo  {journal}
  {Physica D}\ }\textbf {\bibinfo {volume} {238}},\ \bibinfo {pages} {319}
  (\bibinfo {year} {2009})}\BibitemShut {NoStop}%
\bibitem [{\citenamefont {Christlieb}\ \emph {et~al.}(2014)\citenamefont
  {Christlieb}, \citenamefont {Jones}, \citenamefont {Promislow}, \citenamefont
  {Wetton},\ and\ \citenamefont {Willoughby}}]{christlieb2014high}%
  \BibitemOpen
  \bibfield  {author} {\bibinfo {author} {\bibfnamefont {A.}~\bibnamefont
  {Christlieb}}, \bibinfo {author} {\bibfnamefont {J.}~\bibnamefont {Jones}},
  \bibinfo {author} {\bibfnamefont {K.}~\bibnamefont {Promislow}}, \bibinfo
  {author} {\bibfnamefont {B.}~\bibnamefont {Wetton}}, \ and\ \bibinfo {author}
  {\bibfnamefont {M.}~\bibnamefont {Willoughby}},\ }\href@noop {} {\bibfield
  {journal} {\bibinfo  {journal} {Journal of Computational Physics}\ }\textbf
  {\bibinfo {volume} {257}},\ \bibinfo {pages} {193} (\bibinfo {year}
  {2014})}\BibitemShut {NoStop}%
\bibitem [{\citenamefont {Eyre}(1998)}]{eyre1998unconditionally}%
  \BibitemOpen
  \bibfield  {author} {\bibinfo {author} {\bibfnamefont {D.~J.}\ \bibnamefont
  {Eyre}},\ }\href@noop {} {\bibfield  {journal} {\bibinfo  {journal}
  {Unpublished article}\ } (\bibinfo {year} {1998})}\BibitemShut {NoStop}%
\bibitem [{\citenamefont {Doedel}\ \emph {et~al.}()\citenamefont {Doedel},
  \citenamefont {Champneys}, \citenamefont {Fairgrieve}, \citenamefont
  {Kuznetsov}, \citenamefont {Oldeman}, \citenamefont {Paffenroth},
  \citenamefont {Sandstede}, \citenamefont {Wang},\ and\ \citenamefont
  {Zhang}}]{doedel2012auto}%
  \BibitemOpen
  \bibfield  {author} {\bibinfo {author} {\bibfnamefont {E.~J.}\ \bibnamefont
  {Doedel}}, \bibinfo {author} {\bibfnamefont {A.~R.}\ \bibnamefont
  {Champneys}}, \bibinfo {author} {\bibfnamefont {T.}~\bibnamefont
  {Fairgrieve}}, \bibinfo {author} {\bibfnamefont {Y.}~\bibnamefont
  {Kuznetsov}}, \bibinfo {author} {\bibfnamefont {B.}~\bibnamefont {Oldeman}},
  \bibinfo {author} {\bibfnamefont {R.}~\bibnamefont {Paffenroth}}, \bibinfo
  {author} {\bibfnamefont {B.}~\bibnamefont {Sandstede}}, \bibinfo {author}
  {\bibfnamefont {X.}~\bibnamefont {Wang}}, \ and\ \bibinfo {author}
  {\bibfnamefont {C.}~\bibnamefont {Zhang}},\ }\href
  {http://indy.cs.concordia.ca/auto/} {\enquote {\bibinfo {title} {Auto07p:
  Continuation and bifurcation software for ordinary differential equations},}\
  }\BibitemShut {NoStop}%
\bibitem [{\citenamefont {Burke}\ and\ \citenamefont
  {Knobloch}(2007{\natexlab{a}})}]{burke2007homoclinic}%
  \BibitemOpen
  \bibfield  {author} {\bibinfo {author} {\bibfnamefont {J.}~\bibnamefont
  {Burke}}\ and\ \bibinfo {author} {\bibfnamefont {E.}~\bibnamefont
  {Knobloch}},\ }\href@noop {} {\bibfield  {journal} {\bibinfo  {journal}
  {Chaos: An Interdisciplinary Journal of Nonlinear Science}\ }\textbf
  {\bibinfo {volume} {17}},\ \bibinfo {pages} {037102} (\bibinfo {year}
  {2007}{\natexlab{a}})}\BibitemShut {NoStop}%
\bibitem [{\citenamefont {Yochelis}\ \emph
  {et~al.}(2008{\natexlab{a}})\citenamefont {Yochelis}, \citenamefont {Tintut},
  \citenamefont {Demer},\ and\ \citenamefont
  {Garfinkel}}]{yochelis2008formation}%
  \BibitemOpen
  \bibfield  {author} {\bibinfo {author} {\bibfnamefont {A.}~\bibnamefont
  {Yochelis}}, \bibinfo {author} {\bibfnamefont {Y.}~\bibnamefont {Tintut}},
  \bibinfo {author} {\bibfnamefont {L.}~\bibnamefont {Demer}}, \ and\ \bibinfo
  {author} {\bibfnamefont {A.}~\bibnamefont {Garfinkel}},\ }\href@noop {}
  {\bibfield  {journal} {\bibinfo  {journal} {New Journal of Physics}\ }\textbf
  {\bibinfo {volume} {10}},\ \bibinfo {pages} {055002} (\bibinfo {year}
  {2008}{\natexlab{a}})}\BibitemShut {NoStop}%
\bibitem [{\citenamefont {Bergeon}\ \emph {et~al.}(2008)\citenamefont
  {Bergeon}, \citenamefont {Burke}, \citenamefont {Knobloch},\ and\
  \citenamefont {Mercader}}]{PhysRevE.78.046201}%
  \BibitemOpen
  \bibfield  {author} {\bibinfo {author} {\bibfnamefont {A.}~\bibnamefont
  {Bergeon}}, \bibinfo {author} {\bibfnamefont {J.}~\bibnamefont {Burke}},
  \bibinfo {author} {\bibfnamefont {E.}~\bibnamefont {Knobloch}}, \ and\
  \bibinfo {author} {\bibfnamefont {I.}~\bibnamefont {Mercader}},\ }\href@noop
  {} {\bibfield  {journal} {\bibinfo  {journal} {Phys. Rev. E}\ }\textbf
  {\bibinfo {volume} {78}},\ \bibinfo {pages} {046201} (\bibinfo {year}
  {2008})}\BibitemShut {NoStop}%
\bibitem [{\citenamefont {Burke}\ and\ \citenamefont
  {Knobloch}(2007{\natexlab{b}})}]{burke2007snakes}%
  \BibitemOpen
  \bibfield  {author} {\bibinfo {author} {\bibfnamefont {J.}~\bibnamefont
  {Burke}}\ and\ \bibinfo {author} {\bibfnamefont {E.}~\bibnamefont
  {Knobloch}},\ }\href@noop {} {\bibfield  {journal} {\bibinfo  {journal}
  {Physics Letters A}\ }\textbf {\bibinfo {volume} {360}},\ \bibinfo {pages}
  {681} (\bibinfo {year} {2007}{\natexlab{b}})}\BibitemShut {NoStop}%
\bibitem [{\citenamefont {Yochelis}\ \emph
  {et~al.}(2008{\natexlab{b}})\citenamefont {Yochelis}, \citenamefont
  {Knobloch}, \citenamefont {Xie}, \citenamefont {Qu},\ and\ \citenamefont
  {Garfinkel}}]{yochelis2008generation}%
  \BibitemOpen
  \bibfield  {author} {\bibinfo {author} {\bibfnamefont {A.}~\bibnamefont
  {Yochelis}}, \bibinfo {author} {\bibfnamefont {E.}~\bibnamefont {Knobloch}},
  \bibinfo {author} {\bibfnamefont {Y.}~\bibnamefont {Xie}}, \bibinfo {author}
  {\bibfnamefont {Z.}~\bibnamefont {Qu}}, \ and\ \bibinfo {author}
  {\bibfnamefont {A.}~\bibnamefont {Garfinkel}},\ }\href@noop {} {\bibfield
  {journal} {\bibinfo  {journal} {Europhysics Letters}\ }\textbf {\bibinfo
  {volume} {83}},\ \bibinfo {pages} {64005} (\bibinfo {year}
  {2008}{\natexlab{b}})}\BibitemShut {NoStop}%
\bibitem [{\citenamefont {Lloyd}\ \emph {et~al.}(2008)\citenamefont {Lloyd},
  \citenamefont {Sandstede}, \citenamefont {Avitabile},\ and\ \citenamefont
  {Champneys}}]{lloyd2008localized}%
  \BibitemOpen
  \bibfield  {author} {\bibinfo {author} {\bibfnamefont {D.~J.}\ \bibnamefont
  {Lloyd}}, \bibinfo {author} {\bibfnamefont {B.}~\bibnamefont {Sandstede}},
  \bibinfo {author} {\bibfnamefont {D.}~\bibnamefont {Avitabile}}, \ and\
  \bibinfo {author} {\bibfnamefont {A.~R.}\ \bibnamefont {Champneys}},\
  }\href@noop {} {\bibfield  {journal} {\bibinfo  {journal} {SIAM Journal on
  Applied Dynamical Systems}\ }\textbf {\bibinfo {volume} {7}},\ \bibinfo
  {pages} {1049} (\bibinfo {year} {2008})}\BibitemShut {NoStop}%
\bibitem [{\citenamefont {Beaume}\ \emph {et~al.}(2013)\citenamefont {Beaume},
  \citenamefont {Bergeon}, \citenamefont {Kao},\ and\ \citenamefont
  {Knobloch}}]{beaume2013convectons}%
  \BibitemOpen
  \bibfield  {author} {\bibinfo {author} {\bibfnamefont {C.}~\bibnamefont
  {Beaume}}, \bibinfo {author} {\bibfnamefont {A.}~\bibnamefont {Bergeon}},
  \bibinfo {author} {\bibfnamefont {H.-C.}\ \bibnamefont {Kao}}, \ and\
  \bibinfo {author} {\bibfnamefont {E.}~\bibnamefont {Knobloch}},\ }\href@noop
  {} {\bibfield  {journal} {\bibinfo  {journal} {Journal of Fluid Mechanics}\
  }\textbf {\bibinfo {volume} {717}},\ \bibinfo {pages} {417} (\bibinfo {year}
  {2013})}\BibitemShut {NoStop}%
\end{thebibliography}
\end{document}